\documentclass{article}
\usepackage{hyperref}
\usepackage{amssymb}
\usepackage[dvipsnames]{xcolor}
\usepackage[most]{tcolorbox}
\usepackage{amsthm}
\usepackage{amsmath}
\usepackage{physics}
\usepackage{bm}
\usepackage{mathrsfs}
\usepackage{nameref, cleveref}
\usepackage{thmtools}
\usepackage{bbm}
\usepackage[shortlabels]{enumitem}
\usepackage{subcaption}
\usepackage{listings}
\usepackage[page]{appendix}
\usepackage{quoting}
\usepackage[
    style=authoryear,
    maxcitenames=2,
    mincitenames=1,
    uniquelist=false
]{biblatex}

\renewcommand{\P}{\mathbb{P}}

\tcbuselibrary{breakable}

\newtheorem{theorem}{Theorem}[section]
\newtheorem{definition}[theorem]{Definition}

\newtheorem{lemma}[theorem]{Lemma}
\newtheorem{corollary}[theorem]{Corollary}
\newtheorem{assumption}[theorem]{Assumption}
\newtheorem{remark}[theorem]{Remark}

\newcommand{\nCref}[1]{\Cref{#1} (\nameref{#1})}
\newtcolorbox{lemmaborder}{
    blanker,
    borderline west={1.5pt}{0pt}{lightgray},
    left=12pt,
    top=0pt,
    bottom=0pt,
    breakable
}

\crefname{assumption}{assumption}{assumptions}
\crefname{corollary}{Corollary}{Corollaries}
\crefname{remark}{Remark}{Remarks}

\DeclareNameAlias{sortname}{family-given}
\DeclareDelimFormat{multinamedelim}{\addsemicolon\space}
\DeclareDelimFormat{finalnamedelim}{%
     \ifnumgreater{\value{liststop}}{2}{\finalandsemicolon}{}%
     \addspace\bibstring{and}\space}

\title{Optimal Capital Injections in a Microlending Group under External Shocks}
\author{H\'ector Jasso-Fuentes\thanks{Department of Mathematics, Cinvestav - Quer\'etaro, Lib. Norponiente 2000, 76230 Juriquilla, Qro., Mexico. Email: hector.jasso@cinvestav.mx}
\and Alejandra Quintos\thanks{Department of Statistics, University of Wisconsin--Madison, Madison, WI 53706, USA. Email: alejandra.quintos@wisc.edu}
\and Xinta Yang\thanks{Department of Statistics, University of Wisconsin--Madison, Madison, WI 53706, USA. Email: xinta.yang@wisc.edu}}
\date{November 13, 2026}

\hypersetup{
  pdftitle={Optimal Capital Injections in a Microlending Group under External Shocks},
  pdfauthor={H\'ector Jasso-Fuentes, Alejandra Quintos, Xinta Yang}
}

\begin{document}

\maketitle

\begin{abstract}
    Motivated by static models of group default in microlending, we study a dynamic problem in which the members of a lending group are exposed to external shocks, such as droughts or floods, and the manager of a common emergency fund decides at each shock how much capital to inject into the members who were hit. Shocks arrive at random times, the set of members that a shock hits may depend on which members the previous shock hit, and the losses may spill over to members who were not hit directly. Between shocks, the members' wealth follows a linear stochastic differential equation in which members support each other, and the fund grows with interest and member contributions up to a cap. The cost balances the members' shortfall below a poverty line, a welfare measure that contains a standard poverty measure as a special case, against the administrative cost of the injections and the benefit of keeping money in the fund. We formulate the problem as a Markov decision process embedded at the shock times and prove that the embedded process is Markov with a weakly continuous transition kernel. When the discount rate exceeds a growth rate determined by the wealth dynamics, we show that the expected total cost is finite under every policy, that the finite-horizon optimal values converge to the infinite-horizon one at an explicit rate, and that there exists a stationary optimal injection policy, which can be chosen independent of calendar time.
\end{abstract}

\section{Introduction}\label{Sect:Introduction}

Microlending provides small loans to people who have no collateral and no credit history, often through groups whose members are jointly responsible for their loans. The members' businesses are exposed to external shocks, such as droughts, floods or pandemics, which often hit several members at once and may recur: in Bangladesh, a region may be flooded as many as four times in one year \parencite[Chapter~8]{Yunus2007}. Microlending practice has developed reserves at the level of the group against such shocks. In the classic system of the Grameen Bank, every borrower deposited 5\% of each loan in a group fund, from which any member could borrow without interest, with the approval of all the other members, for instance to rebuild a business hit by a natural disaster \parencite[Chapter~5]{Yunus2007}. This raises the question that we study: when a shock hits some members of a group, how much should be injected into them from a common emergency fund, given that the money spent now is no longer available for the next shock?

Mathematical models of microlending have studied access to credit, loan pricing and default. \textcite{DienerDienerKhodrProtter2009} and \textcite{BoiquayeKorvehGyamerah2026} describe the credit status of a borrower or of a group by discrete-time Markov chains, \textcite{JarrowProtter2019} determine fair interest rates for group loans, and \textcite{BoiquayeProtter2024} study the probability of no default when the members' wealth follows independent geometric Brownian motions. \textcite{ProtterQuintos2022} and \textcite{JassoFuentesQuintosYang2026} study the probability that a group defaults over a single financing period, without and with contagion among members. None of these models follows the members' wealth over time together with external shocks and a common fund. Our model is motivated by the last two papers, but it is a different, dynamic model, in which a common fund is used over time to support the members who are hit by shocks.

In our model (Section~\ref{Sect:Problem.Setup}), the members' wealth follows, between shocks, a linear stochastic differential equation in which the members support each other. Shocks arrive at random times, and each shock directly hits a set of members, its mark. The time and the mark of the next shock depend on the past only through the most recent mark, so that, for instance, a drought may follow a drought. At a shock, members lose a random fraction of their positive wealth, and the losses may spill over to members who were not hit directly. The group has an emergency fund, which grows with interest and with the members' contributions up to a cap. At each shock, the manager of the fund, either the group itself or the microfinance institution, decides how much to inject into the members who were directly hit, without exceeding the fund's balance. The cost at each shock has three terms: a welfare term, which penalizes each member's shortfall below a poverty line through a function $f$; the administrative cost of the injections; and a reward for keeping money in the fund. For instance, if all members have the same poverty line $l > 0$ and the same weight $1/(dl)$, and $f(z) = \max(z, 0)$, so that the penalty is the shortfall itself, then the welfare term is the poverty measure $P_1$ of \textcite{FosterGreerThorbecke1984}, applied to the members' wealth after the injections.

We formulate the problem as a Markov decision process embedded at the shock times (Section~\ref{Sect:Existence.Optimal.Policy}). We prove that, under every policy, the embedded process is Markov, with a weakly continuous transition kernel that does not depend on the policy (\Cref{theorem:markov-property,theorem:transition-kernel,theorem:weak-continuity}). The cost is unbounded and takes both signs. When the discount rate exceeds a growth rate determined by the wealth dynamics, the expected total cost is finite under every policy, and the optimal costs over $N$ shocks converge to the optimal total cost at the explicit rate $\mathscr{O}(Nq^N)$, with $q \in (0, 1)$ (\Cref{theorem:one-stage-cost-estimate,corollary:V-defined,corollary:V_N-convergence-rate}). Finally, there exists a stationary optimal injection policy, and it can be chosen independent of calendar time (\Cref{theorem:existence-of-optimal-policy,remark:time-independence}).

Mathematically, our problem belongs to a class of control problems in which the controller may act only when an exogenous signal arrives. Examples are the explicit solution of a stopping problem for a geometric Brownian motion by \textcite{DupuisWang2002}; optimal stopping and impulse control of Markov--Feller processes by Menaldi and Robin \parencite*{MenaldiRobin2016,MenaldiRobin2017}; and optimal stopping of general Markov processes, through an equivalent Markov decision process, by \textcite{JassoFuentesVasquez2025}; see also the references in these papers. In these works, a signal only opens an opportunity to act, and the general results are proved for bounded costs or rewards. Here, the signals are the shocks themselves: they change the state, and their marks, which depend on the previous mark, determine which members are eligible for injections. Moreover, the injections are limited by a fund with its own dynamics, and the cost is unbounded and takes both signs. Like \textcite{JassoFuentesVasquez2025} and \textcite{AlbrecherBauerleThonhauser2011}, we work with the process observed at the signal times; this approach is also standard for piecewise deterministic Markov processes \parencite{Davis1993}, in which, unlike in our model, the dynamics between jumps are controlled. The features of our problem listed above require additional arguments. Since each shock changes the state, we obtain the Markov property and the weak continuity of the kernel from a representation of the wealth between shocks by a single functional that does not depend on the stage or on the policy \parencite{Kallenberg}. Since the cost is unbounded, we bound it uniformly over policies by combining the growth rate of the wealth with the Laplace transform of the times between shocks; a bound of the same kind, for a one-dimensional surplus, is in \textcite[Lemma~3.1]{AlbrecherBauerleThonhauser2011}. Since the cost takes both signs, we subtract a term that does not depend on the policy, which makes the cost nonnegative, so that the existence results of \textcite{Feinberg} apply.

Related problems appear in actuarial science, where dividends or capital injections are restricted to random times independent of the surplus, often the arrival times of a Poisson process \parencite{AlbrecherCheungThonhauser2011,AlbrecherBauerleThonhauser2011,AvramPerezYamazaki2018,NobaPerezYamazakiYano2018,MataMorenoFrancoNobaPerez2022}; there the surplus is one-dimensional, the decision times are independent of the surplus and of its jumps, and the results are explicit identities or optimal strategies of barrier or band type.
Our peer-support term has the same form as the coupling through lending in models of interbank lending \parencite{FouqueIchiba2013,CarmonaFouqueSun2015}, in which there are no shocks and no common fund that is used at the shocks.

This paper sets up the model and shows that the optimization problem is well posed; numerical results will appear in a later version. Section~\ref{Sect:Problem.Setup} presents the model, Section~\ref{Sect:Existence.Optimal.Policy} the embedded Markov decision process and our results, and the appendix collects auxiliary results.


\section{Problem Setup}\label{Sect:Problem.Setup}

We consider a lending group of $d$ borrowers who receive microloans under group lending with joint liability.
Each member runs a small business, and several members may operate in the same sector.
We denote by $Y^j_t$ the business wealth of member $j$ at time $t$; a negative value of $Y^j_t$ represents outstanding debt.
Between shocks, the wealth vector $Y_t = (Y^1_t, \dots, Y^d_t)$ follows an SDE with drift coefficient function $b(\cdot)$ and diffusion coefficient function $\sigma(\cdot)$, driven by a $d$-dimensional Brownian motion $B$ that captures the day-to-day fluctuations of the members' businesses; the specific coefficients, which include support transfers among members, are given in \cref{equation:wealth-dynamics-1} below.
In addition, the members' businesses are exposed to external shocks to their sectors, such as a drought, a pest outbreak, or a collapse in local prices, which cause sudden losses in the wealth of the members affected.
To protect the group, capital from an emergency fund, built from member contributions and interest, can be injected into the members who were hit by a shock, in the form of emergency loans or grants.
Our goal is to determine how much to inject at each shock, balancing the welfare of the members who were hit against the need to preserve the fund for future shocks.

We assume that the Brownian motion $B$ is independent of the sector shocks, that is, of their timing, of which members they hit, and of their severity.
Sector shocks therefore affect the members' wealth only through the jumps described in the next subsection.
This is reasonable when each member's business is small relative to its sector, as is typical of microfinance borrowers (for instance, small farms or shops in rural villages): the ups and downs of a single business have no influence on shocks at the sector level, and the effect of a sector shock on the business is captured by the jump in wealth at the shock time.

\emph{Notation.} For $y \in \mathbb{R}^d$, $||y||$ denotes the Euclidean norm and $||y||_1 := \sum_{j=1}^d |y_j|$; for a matrix $M = (M_{i,j})$, $||M|| := \big(\sum_{i,j} M_{i,j}^2\big)^{1/2}$ denotes the Frobenius norm. We write $\mathbb{R}_+ := [0, \infty)$, $z^+ := \max(z, 0)$, $z \wedge z' := \min(z, z')$ and $z \vee z' := \max(z, z')$ for $z, z' \in \mathbb{R}$.

\subsection{Dynamics Jumping Mechanism}

In this subsection, we describe when shocks occur, which members they hit, and how the members' wealth jumps at those moments. We call the arrival of a shock a \emph{market signal}, and the set of members it directly hits the \emph{mark} of the signal.
Marks take values in
\[
    \mathscr{K} := \{k : k \subseteq \{1, \dots, d\},\ k \neq \emptyset\},
\]
a finite set with $2^d - 1$ elements, which we equip with the discrete topology. For example, the mark $\{3\}$ means that only member three is directly hit, and the mark $\{1, 2\}$ means that members one and two are directly hit simultaneously.

    The timing and the marks of the market signals are governed by nonnegative rates $\lambda_{k,k'}$, for $k, k' \in \mathscr{K}$. For each $k \in \mathscr{K}$, let
    \begin{equation} \label{equation:signal-total-rate}
        \lambda_k := \sum_{k' \in \mathscr{K}} \lambda_{k,k'}, \qquad \bar{\lambda} := \max_{k \in \mathscr{K}} \lambda_k,
    \end{equation}
    and assume that $\lambda_k > 0$ for every $k \in \mathscr{K}$.
    We assume that the problem starts right after a market signal, whose mark we denote by $k_0 \in \mathscr{K}$.
    Let $(E_n)_{n \ge 1}$ be i.i.d.\ exponential random variables with mean $1$, and let $(U_n)_{n \ge 1}$ be i.i.d.\ random variables uniformly distributed on $[0, 1)$.
    Fix an enumeration $k^1, \dots, k^{2^d-1}$ of $\mathscr{K}$, and define $F: \mathscr{K} \times [0, 1) \rightarrow \mathscr{K}$ by $F(k, u) := k^m$ whenever $\sum_{i < m} \lambda_{k,k^i}/\lambda_k \le u < \sum_{i \le m} \lambda_{k,k^i}/\lambda_k$.
    Set $\tau_0 := 0$, $K_0 := k_0$ and, for $n \ge 0$,
    \begin{equation} \label{equation:signal-construction}
        \tau_{n+1} := \tau_n + \frac{E_{n+1}}{\lambda_{K_n}}, \qquad K_{n+1} := F(K_n, U_{n+1}).
    \end{equation}
    For $n \ge 1$, $\tau_n$ is the time of the $n$-th market signal and $K_n$ is its mark.
    Finally, let $\kappa_t := K_n$ for $\tau_n \le t < \tau_{n+1}$, so that $\kappa_t$ is the mark of the most recent market signal at time $t$ and $\kappa_{\tau_n} = K_n$.
    Since $\lambda_{K_n} \le \bar{\lambda}$, we have $\tau_n \ge \bar{\lambda}^{-1} \sum_{i=1}^n E_i \rightarrow \infty$ a.s., so only finitely many signals occur in any bounded time interval.

    In words: given that the most recent signal hit the set of members $k$, the time until the next signal is exponential with rate $\lambda_k$, and, independently of this time, the next signal hits the set $k'$ with probability $\lambda_{k,k'}/\lambda_k$ (see \Cref{lemma:more-on-signal-interval-distribution}). This is the jump chain and holding time construction of a continuous-time Markov chain in \textcite[Section~2.6]{Norris}, except that the next mark may coincide with the current one. Equivalently, after a signal with mark $k$, each set $k'$ has an independent exponential clock with rate $\lambda_{k,k'}$, and the next signal hits the set whose clock rings first; see \textcite[Theorem~2.3.3]{Norris}.
    In particular, since $\tau_{n+1} - \tau_n$ and $K_{n+1}$ are built only from $K_n$, $E_{n+1}$ and $U_{n+1}$, the time until the next signal and the set of members it hits depend on the past only through the mark of the most recent signal; they do not depend on the members' wealth, the fund's balance, or time.

The rates allow the next shock to depend on the previous one: after a signal that hits the set $k$, the next signal hits the set $k'$ with probability $\lambda_{k,k'}/\lambda_k$. For instance, consider a group with two farmers (members $1$ and $2$) and a herder (member $3$). A drought may hit both farmers, $\{1, 2\}$, and be followed by losses for the herder, $\{3\}$, as fodder becomes scarce; this is captured by a large $\lambda_{\{1,2\},\{3\}}$. A drought may also be followed by another drought, which corresponds to $\lambda_{\{1,2\},\{1,2\}} > 0$, and a flood may hit every member, $\{1, 2, 3\}$. If the two farmers are always hit together, that is, every shock that hits one of them also hits the other, then the marks $\{1\}$, $\{2\}$, $\{1,3\}$ and $\{2,3\}$ never occur, which corresponds to $\lambda_{k,k'} = 0$ for these $k'$ and every $k \in \mathscr{K}$.

\begin{remark} \label{remark:poisson-representation}
    When $\lambda_{k,k} = 0$ for all $k \in \mathscr{K}$, the process $\kappa$ is a continuous-time Markov chain on $\mathscr{K}$ with jump rates $\lambda_{k,k'}$, and the $\tau_n$, $n \ge 1$, are its jump times.
    After identifying $\mathscr{K}$ with $\{1, \dots, 2^d - 1\}$ through an enumeration, such a chain can equivalently be written as $\dd{\kappa_t} = \int_{\mathbb{R}} h(\kappa_{t-}, z)\, \mathfrak{p}(\dd{t}, \dd{z})$ for a Poisson random measure $\mathfrak{p}$ on $\mathbb{R}_+ \times \mathbb{R}$ with intensity $\dd{t} \times \dd{z}$ and a suitable function $h$; see \textcite[Chapter~II, \S 2.1, eq.~(41)]{Skorohod} and \textcite[eqs.~(3.3) and (3.4)]{Ghosh}.
    This representation is common in the literature on switching diffusions.
    When $\lambda_{k,k} > 0$ for some $k \in \mathscr{K}$, a signal with $K_{n+1} = K_n$ does not change $\kappa$, so the signal times are no longer the jump times of $\kappa$; this is why we construct the signals directly. (The process $\kappa$ is still a continuous-time Markov chain, with jump rates $\lambda_{k,k'}$ for $k' \neq k$, but it does not record a signal that hits the same set of members as the previous one.)
\end{remark}

At a signal time, the members who were directly hit, and possibly other members, suffer a sudden loss in wealth. We use a collection of i.i.d. random matrices $(\eta_n)_{n \ge 1}$, each taking values in $[0, 1)^{|\mathscr{K}|\times d}$, to describe the fractions of wealth that the members lose at time $\tau_n$. The rows of $\eta_n$ are indexed by the marks $k \in \mathscr{K}$ and its columns by the members $j \in \{1, \dots, d\}$, and we denote its entries by $\eta_{n,k,j}$: the entry $\eta_{n,k,j}$ is the fraction of wealth that member $j$ loses at time $\tau_n$ if the signal has mark $k$. Only the row of the mark that actually occurs, $k = K_n$, is used. We allow one row per mark because the losses may depend on which members were hit and may differ across members.
For example, under the mark $\{1, 2\}$, members $1$ and $2$ would typically lose a larger fraction of their wealth than the other members, which corresponds to larger entries $\eta_{n,k,j}$ for $k = \{1, 2\}$ and $j \in \{1, 2\}$.
To be precise, we shall denote by $G: \mathbb{R}^d \times \mathscr{K} \times [0, 1)^{|\mathscr{K}|\times d} \rightarrow \mathbb{R}^d$ the shock effect function, which is defined via
\begin{equation} \label{equation:jump-down-mechanism-1}
    G(y, k, u) := \big(y_1 - u_{k,1}\,y_1^+, \dots, y_d - u_{k,d}\,y_d^+\big).
\end{equation}
Here $u = (u_{k',j})_{k' \in \mathscr{K},\, 1 \le j \le d} \in [0, 1)^{|\mathscr{K}| \times d}$ is a matrix with rows indexed by the marks and columns indexed by the members, just as for $\eta_n$, and $y_j^+ := \max(y_j, 0)$. Thus the $j$-th component of $G(y, k, u)$ uses only the row of $u$ that corresponds to the mark $k$, and at a signal time $u$ is replaced by $\eta_n$ and $k$ by $K_n$.

Notice $G$ is a 1-Lipschitz function and is concave in variable $y$.
Also, under $G$ a member who is already in debt, i.e. whose current wealth is negative, does \textit{not} see the debt increase at a shock; a member with positive wealth loses the fraction $u_{k,j}$ of it. Note that $G$ allows $u_{k,j} > 0$ also for $j \notin k$: members who were not directly hit may still lose wealth through spillovers within the group, in the spirit of the contagion studied in \textcite{JassoFuentesQuintosYang2026}. We call such members \emph{indirectly affected}. If only the members directly hit should lose wealth, it suffices to assume that $\eta_{n,k,j} = 0$ for $j \notin k$. Immediately after the shock, capital may be injected from the fund, as emergency loans or grants, \textit{only} into the members who were directly hit, i.e. the members $j \in \kappa_{\tau_n}$; indirectly affected members are not eligible.

Next, we introduce the group's emergency fund, whose balance at time $t$ we denote by $R_t$. The fund may be held by the group itself or by the microfinance institution on the group's behalf. We call whoever manages the fund the \emph{fund's manager}; the fund's manager decides how much to inject at each shock. Between market signals it accumulates without risk, according to
\begin{equation} \label{equation:money-reserve-1}
    dR_t = \left( \rho R_t + \sum_{j=1}^d \delta_j \max(Y_t^j, 0) \right) \cdot 1_{\{R_t < \bar{r}\}} dt,
\end{equation}
where $\rho \ge 0$ is the interest rate earned on the fund's balance, $\delta_j \ge 0$ is the rate at which member $j$ contributes to the fund while their wealth is positive (members in debt do not contribute), and $\bar{r}$ is the cap of the fund. At time $0$, before the injection $\xi_0$, the fund holds a given amount $0 \le r_0 \le \bar{r}$. The cap $\bar{r}$ is the target size of the fund. Once the balance reaches $\bar{r}$, the members stop contributing and interest is no longer added to the fund (it may, for instance, be returned to the members or used by the lender elsewhere). This keeps the fund from growing indefinitely at the members' expense: money set aside beyond what is needed for emergencies is money the members cannot invest in their businesses.

At each signal time $\tau_n$, $n \ge 0$, the injection $\xi_n$ is paid out of the fund, so that $R_{\tau_n} = R_{\tau_n-} - ||\xi_n||_1$, with $R_{0-} := r_0$; between signals, $R$ evolves according to \cref{equation:money-reserve-1}, starting from $R_{\tau_n}$. Since an injection cannot exceed the fund's balance (\Cref{definition:feasible-action-set}), the balance never becomes negative. Moreover, if $0 \le R_{\tau_n} \le \bar{r}$, then \cref{equation:money-reserve-1} has a unique solution on $[\tau_n, \tau_{n+1})$, namely $\mathfrak{R} \wedge \bar{r}$, where $\mathfrak{R}$ solves the same equation without the indicator, started from $R_{\tau_n}$ (see \Cref{lemma:ODE-solution}). Indeed, $\mathfrak{R}$ is nondecreasing because $\rho \ge 0$ and $\delta_j \ge 0$, so $\mathfrak{R} \wedge \bar{r}$ coincides with $\mathfrak{R}$ until the balance reaches $\bar{r}$ and then stays at $\bar{r}$, where the right-hand side of \cref{equation:money-reserve-1} vanishes; conversely, any solution coincides with $\mathfrak{R}$ while it is below $\bar{r}$ and cannot increase once it reaches $\bar{r}$. In particular, $0 \le R_t \le \bar{r}$ for all $t \ge 0$ (see the proof of \Cref{lemma:functional-rep-of-R}).

Finally, we specify the members' wealth dynamics $Y$ between consecutive market signals. We use a linear SDE with multiplicative noise, in which the members are coupled through support transfers:
\begin{equation} \label{equation:wealth-dynamics-1}
    \dd{Y_t^j} = \left( \mu_j Y_t^j + \sum_{l \neq j} \theta_{j,l} (Y_t^l - Y_t^j) \right) \dd{t} + \sigma_j Y_t^j \dd{B_t^j},
\end{equation}
where $\mu_j$ is the growth rate of member $j$'s business, $\theta_{j,l} \ge 0$ measures the support that member $j$ receives from member $l$, and $\sigma_j$ is the volatility of member $j$'s business. Without the coupling, that is, if $\theta_{j,l} = 0$ for all $j \ne l$, each $Y^j$ would be a geometric Brownian motion. The term $\theta_{j,l}(Y^l_t - Y^j_t)$ pulls member $j$'s wealth toward member $l$'s, so a member who falls behind receives support from better-off members. It models the mutual assistance among members that group lending with joint liability encourages. The same coefficients also measure how much a member is exposed to the others' losses: the term $\theta_{l,j}(Y^j_t - Y^l_t)$ in member $l$'s equation is negative when member $j$ is worse off. Hence a loss suffered by one member, for instance at a shock that hits only that member, gradually lowers the wealth of the members who support them, and an injection into a member who was hit also benefits, indirectly, the members who support them, although these members are not eligible for injections themselves. If all members start with nonnegative wealth, their wealth stays nonnegative at all times: the losses at a shock are proportional to positive wealth, injections are nonnegative, and at a member with zero wealth the noise vanishes and the drift $\sum_{l \ne j} \theta_{j,l} Y^l_t$ is nonnegative. Debt in the model therefore comes only from the initial condition (for instance, an outstanding balance at time $0$), and peer support can then spread it to other members.
We also provide its matrix form.
Under the notation
\begin{equation} \label{equation:wealth-dynamics-2}
    \Theta_{j,j} := \mu_j - \sum_{l \ne j} \theta_{j,l}, \quad \Theta_{j,l} := \theta_{j,l} \text{ for } j \ne l,
\end{equation}
the wealth dynamics between adjacent market signals is characterized by
\begin{equation} \label{equation:wealth-dynamics-3}
    \dd{Y_t} = \Theta Y_t\dd{t} + \text{diag}(\sigma_1Y_t^1, ..., \sigma_dY_t^d) \dd{B_t}.
\end{equation}
Thus the drift and diffusion coefficients $b$ and $\sigma$ mentioned at the beginning of this section are $b(y) := \Theta y$ and $\sigma(y) := \text{diag}(\sigma_1 y_1, \dots, \sigma_d y_d)$ for $y \in \mathbb{R}^d$; both are Lipschitz continuous and have linear growth.

    All random variables are defined on a complete probability space $(\Omega, \mathcal{F}, \P)$ that supports the Brownian motion $B$, the sequences $(E_n)_{n \ge 1}$, $(U_n)_{n \ge 1}$ and $(\eta_n)_{n \ge 1}$ introduced above, and a sequence $(\chi_n)_{n \ge 0}$ of i.i.d.\ random variables uniformly distributed on $[0, 1)$, which the fund's manager may use to randomize the injections (see \Cref{definition:policies}); all of these are independent (\Cref{assumption:standing}).
    Let $\mathcal{G} := \sigma(E_n, U_n, \eta_n, \chi_m : n \ge 1, m \ge 0)$, and let $(\mathcal{F}_t)_{t \ge 0}$ be the usual augmentation of $(\sigma(B_s : s \le t) \vee \mathcal{G})_{t \ge 0}$.
    Since $B$ is independent of $\mathcal{G}$, it is an $(\mathcal{F}_t)$-Brownian motion.
    Each $\tau_n$ is $\mathcal{G}$-measurable by \cref{equation:signal-construction}, hence $\mathcal{F}_0$-measurable and, in particular, an $(\mathcal{F}_t)$-stopping time.
    The filtration $(\mathcal{F}_t)$ contains information about future signals; it is only a device to construct the wealth process. What the fund's manager knows when choosing an injection is described by the class of policies in \Cref{definition:policies}.

We now provide a piecewise representation of process $Y$.
First, define the time-shifted filtration and Brownian motion via
\begin{equation*}
    \mathbb{F}^{(n)} := (\mathcal{F}^{(n)}_t)_{t \ge 0}, \quad \mathcal{F}^{(n)}_t := \mathcal{F}_{\tau_n+t}
\end{equation*}
and
\begin{equation*}
        B^{(n)}_\cdot := B_{\tau_n + \cdot} - B_{\tau_n}.
\end{equation*}
By the strong Markov property, $B^{(n)}$ is an $\mathbb{F}^{(n)}$-Brownian motion independent of $\mathcal{F}_{\tau_n}$.
Next, under filtration $\mathbb{F}^{(n)}$, consider the SDE
\begin{equation*}
    \dd{Y^{(n)}_t} = \Theta Y^{(n)}_t \dd{t} + \text{diag}\left(\sigma_1Y^{(n), 1}_t, ..., \sigma_dY^{(n), d}_t \right) \dd{B^{(n)}_t},
\end{equation*}
where matrix $\Theta$ is defined in \cref{equation:wealth-dynamics-2}.
Let $y_0 \in \mathbb{R}^d$ be the members' wealth at time $0$, right after the market signal with mark $k_0$, so that the losses of that signal are already included in $y_0$, and let $\xi_n$ be the injection at time $\tau_n$, chosen by a policy (\Cref{definition:policies}). Process $Y^{(0)}$ starts at $Y^{(0)}_0 := y_0 + \xi_0$, and for $n \ge 1$ we define inductively $Y^{(n)}_0 := G(Y^{(n-1)}_{\tau_n-\tau_{n-1}}, \kappa_{\tau_n}, \eta_n)+\xi_n$, which is $\mathcal{F}_{\tau_n}$-measurable.
Finally, define
\begin{equation} \label{equation:Y-definition}
    Y_t := \sum_{n=0}^\infty 1_{[\tau_n, \tau_{n+1})}(t) Y^{(n)}_{(t-\tau_n)\vee 0}.
\end{equation}
Because $\tau_n \nearrow \infty$ a.s., for any given time $t$, the summation above only has finitely many nonzero terms, so it is well defined.

\begin{remark} \label{remark:global-formula-of-Y}
    Process $Y$ defined above satisfies
    \begin{multline} \label{equation:theorem-global-formula-of-Y-1}
            Y_t = y_0 + \xi_0 + \int_0^t \Theta Y_{s-} \dd{s} + \int_0^t \text{diag}(\sigma_1 Y^1_{s-}, ..., \sigma_d Y^d_{s-}) \dd{B_s} \\
            + \sum_{n \ge 1,\, \tau_n \le t}  [G(Y_{\tau_n-}, \kappa_{\tau_n}, \eta_n) + \xi_n - Y_{\tau_n-}],
    \end{multline}
    in the sense that $Y$ is indistinguishable from the process given on the right hand side above.
    The right-hand side of \cref{equation:theorem-global-formula-of-Y-1} has the same structure as the controlled jump diffusion in equation (1.2) of \textcite{Menaldi1}: the jumps $G(Y_{\tau_n-}, \kappa_{\tau_n}, \eta_n) - Y_{\tau_n-}$ come from the market signals, and the injections $\xi_n$ play the role of the controller's impulses. Unlike in \textcite{Menaldi1}, where the controller also chooses when to apply the impulses, here injections can only be made at the signal times.
    Throughout the rest of this paper though, we will follow the piecewise formulation of $Y$.

    The idea of representing a global process piecewise can be found in the proof of \textcite[Theorem~3.13]{Mao}, or at the end of Theorem V.2.4's proof in \textcite{Ikeda}.
    The following diagram illustrates how $Y$ evolves over time.
\end{remark}

\begin{figure}[h]
    \centering
    \includegraphics[width=0.9\linewidth]{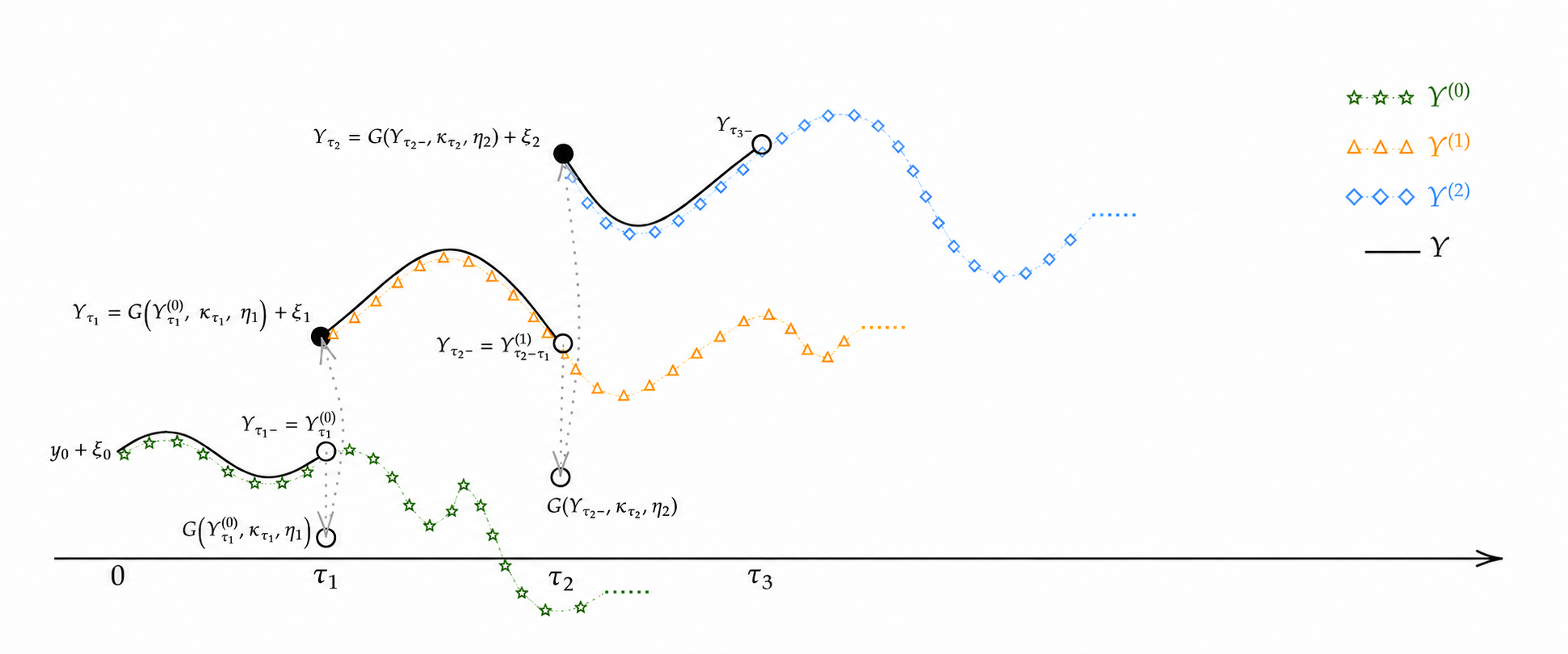}
    \caption{The members' wealth $Y$ over time (downward gray dotted arrows indicate the losses at a market signal, and upward arrows indicate injections from the emergency fund)}
    \label{fig:Y-dynamics}
\end{figure}

Before we define the cost, we introduce the state space of the system and the set of feasible injections.
\begin{definition}[System State Space] \label{definition:system-state-space}
    The system's state value consists of four components: the members' current wealth, the mark of the most recent market signal (the set of members directly hit), the fund's balance, and time.
    Therefore, define $\mathbb{X} := \mathbb{R}^d \times \mathscr{K} \times [0, \bar{r}] \times \mathbb{R}_+$ to be the system's state space.
    Notation $x = (y, k, r, t) \in \mathbb{X}$ shall appear frequently throughout the rest of this paper.
\end{definition}

\begin{definition}[Feasible Action Set] \label{definition:feasible-action-set}
    Given the system's value $x = (y, k, r, t)$, define the feasible action set to be
    \begin{equation} \label{equation:definition-feasible-action-set-1}
        A(x) := \{a \in \mathbb{R}_+^d : ||a||_1 \le r, a_j = 0 \text{ for } j \notin k\}.
    \end{equation}
    It enforces two things: the total injection never exceeds the fund's balance, and only the members directly hit by the market signal are eligible for an injection. The graph of $A$ is
    \begin{equation} \label{equation:graph-of-A}
        \text{Gr}(A) := \{(x, a) : x \in \mathbb{X},\ a \in A(x)\}.
    \end{equation}
    Since $\mathscr{K}$ is discrete, $\text{Gr}(A)$ is a closed subset of $\mathbb{X} \times \mathbb{R}_+^d$.
\end{definition}

The task now is to define the cost function $c(x, a)$, which leads to our optimization target defined in the next section.
\begin{equation} \label{equation:cost-function}
    c(x, a) := e^{-\gamma t}\left[\sum_{j=1}^d w_jf(l_j-y_j-a_j) + w_A||a||_1 - w_R(1+r-||a||_1)^\nu\right],
\end{equation}
where $\gamma > 0$ is the discount rate, $w_A, w_R, w_j \ge 0$ are given weight constants, $f : \mathbb{R} \rightarrow [0, \infty)$ is a welfare penalty, continuous and with at most linear growth, and the exponent $\nu \in (0, 1)$; the precise conditions on these constants and on $f$ are in \Cref{assumption:standing} below.
In terms of variables, $x=(y, k, r, t) \in \mathbb{X}$ is the system's state (\Cref{definition:system-state-space}) and $a \in \mathbb{R}_+^d$ is the action, i.e. the vector of injections into the members. We evaluate $c$ only at pairs $(x, a)$ with $||a||_1 \le r$, in particular at pairs in the graph $\text{Gr}(A)$ of the feasible action set (\Cref{definition:feasible-action-set}); then $1 + r - ||a||_1 \ge 1$, so the power in the third term is well defined.
The motivation for presenting such a cost function is as follows.
We start with a discounting factor depending on time.
Inside the square bracket,
the first term measures the members' welfare shortfall. The constant $l_j$ is a poverty line, or subsistence wealth level, for member $j$, and the term penalizes post-injection wealth $y_j + a_j$ below it.
The argument $l_j - y_j - a_j$ of $f$ is the shortfall of member $j$'s post-injection wealth below the poverty line. It is natural to take $f$ nondecreasing, so that a larger shortfall is penalized more, although none of our results requires this. A standard choice is the positive part $f(z) = z^+ := \max(z, 0)$, for which the welfare term is a weighted sum of the members' shortfalls below the poverty line, as in the poverty measure $P_1$ of \textcite{FosterGreerThorbecke1984}: once an injection lifts a member's wealth to the poverty line, the penalty for that member vanishes. Another choice is the smooth approximation of the positive part given by $f(z) = f_\beta(z) := \frac{1}{\beta}\log(1 + e^{\beta z})$ with $\beta > 0$, which converges to $z^+$ as $\beta \rightarrow \infty$ and may be convenient in numerical computations. Both satisfy the conditions in \Cref{assumption:standing}, with $C_f = 1$ and $C_f = \max(\log 2/\beta, 1)$, respectively.
The second term inside the square bracket is the administrative cost of disbursing emergency loans or grants, proportional to the total amount injected.
The third term rewards keeping money in the fund for future shocks: a larger remaining balance $r-||a||_1$ leads to a smaller total cost $c(x, a)$, due to the leading minus sign.
This third term is motivated by the classic constant relative risk aversion (CRRA) utility function.
The reason to pick the power $\nu \in (0, 1)$ is that the fund has a maximum size $\bar{r}$, and the concavity of $(\cdot)^\nu$ captures the idea that, as the remaining balance gets closer to this maximum, the reward grows more and more slowly.
Finally, the $w_j, w_A, w_R$ put weights on those quantities since we might emphasize different components of the total cost.
To summarize, the fund's manager faces a trade-off between supporting the members who were hit now and preserving the fund's balance for future shocks.

We close this section by listing our model's assumptions.

\begin{assumption}[Standing assumptions] \label{assumption:standing}
    \begin{enumerate}[(i)]
        \item $\lambda_k > 0$ for every $k \in \mathscr{K}$.
         \item $\mu_j \in \mathbb{R}$, $\sigma_j \in \mathbb{R}$ and $\theta_{j,l} \ge 0$ for all $j \ne l$ are constants.
        \item $\rho \ge 0$, $\delta_j \ge 0$ for every $j$, $\bar{r} > 0$, and $0 \le r_0 \le \bar{r}$.
        \item $\gamma > 0$, $\nu \in (0, 1)$, $w_j \ge 0$ and $l_j \ge 0$ for every $j$, and $w_A, w_R \ge 0$.
        \item The welfare penalty $f$ in \cref{equation:cost-function} is continuous, takes values in $[0, \infty)$, and grows at most linearly: there is a constant $C_f \ge 0$ such that $f(z) \le C_f(1 + z^+)$ for all $z \in \mathbb{R}$.
        \item The Brownian motion $B$, the sequences $(E_n)_{n \ge 1}$, $(U_n)_{n \ge 1}$ and $(\chi_n)_{n \ge 0}$, and the i.i.d.\ shock matrices $(\eta_n)_{n \ge 1}$ are independent.
    \end{enumerate}
\end{assumption}

\section{Existence of Optimal Policy}\label{Sect:Existence.Optimal.Policy}
In this section, we prove the existence of a stationary deterministic policy $\varpi^*$ that satisfies $V(\varpi^*, x) = V^*(x)$ for all $x \in \mathbb{X}$.
We start this section with a direct application of Theorem 1 in \textcite{Kallenberg}.

\begin{lemma}[Functional Solution to SDE] \label{lemma:functional-solution-to-SDE}
    For $Y^{(n)}$ and $B^{(n)}$ defined preceding \cref{equation:Y-definition},
    there exists a Borel measurable functional $\mathfrak{F}:\mathbb{R}^d \times \mathscr{C}(\mathbb{R}_+, \mathbb{R}^d) \rightarrow \mathscr{C}(\mathbb{R}_+, \mathbb{R}^d)$ satisfying
    \begin{equation} \label{equation:lemma-functional-solution-to-SDE-1}
        Y^{(n)} = \mathfrak{F}(Y^{(n)}_0, B^{(n)}), \quad \forall n \in \mathbb{N}
    \end{equation}
    in the sense of indistinguishability.
\end{lemma}
\begin{proof}
    Since $b$ and $\sigma$ are Lipschitz continuous, \cref{equation:wealth-dynamics-3} has a strong solution for every deterministic initial condition, and pathwise uniqueness holds; in particular, weak existence and pathwise uniqueness hold for solutions starting at fixed points. Theorem 1 of \textcite{Kallenberg} therefore gives a Borel measurable functional $\mathfrak{F}$ such that every solution of \cref{equation:wealth-dynamics-3}, on any filtered probability space and with any initial distribution, is almost surely equal to $\mathfrak{F}$ applied to its initial value and its driving Brownian motion. For each $n$, $B^{(n)}$ is an $\mathbb{F}^{(n)}$-Brownian motion, and the equation for $Y^{(n)}$ is \cref{equation:wealth-dynamics-3} with $B^{(n)}$ in place of $B$, so $Y^{(n)} = \mathfrak{F}(Y^{(n)}_0, B^{(n)})$ a.s.; since both sides have continuous paths, they are indistinguishable. The functional $\mathfrak{F}$ depends neither on $n$ nor on the distribution of $Y^{(n)}_0$, which depends on $n$ and on the policy. This is why we use Theorem 1 of \textcite{Kallenberg} rather than the classical Yamada--Watanabe theorem, in which the functional may depend on the initial distribution (see the introduction of \textcite{Kallenberg}).
\end{proof}

\begin{corollary}
    Process $Y^{(n)}$ satisfies $Y^{(n)}_t \in \sigma(Y^{(n)}_0)\vee \sigma(B^{(n)}_s; s \ge 0)$, where ``$\vee$" means the sigma-field generated by the two sets.
\end{corollary}

Intuitively, this corollary states that given the initial condition $Y^{(n)}_0$, the rest of $Y^{(n)}$ is independent of the history before time $\tau_n$ (since $B^{(n)}$ is), which hints the Markov property that we shall closely examine later.
The first task now is to discretize the controlled dynamics $Y$ and the fund's balance $R$ at signal times.

\begin{definition}[Discretized Dynamics] \label{definition:discretized-dynamics}
    For any given $(y_0, k_0, r_0, t_0) \in \mathbb{X}$, define
    \begin{equation} \label{equation:definition-discretized-dynamics-1}
        X_0 := (y_0, k_0, r_0, t_0)
    \end{equation}
    to be the initial condition. Here $t_0$ is the calendar time at which the problem starts: the processes $Y$ and $R$ and the signal times $\tau_n$ are indexed by the time elapsed since $t_0$, so the $n$-th signal occurs at calendar time $t_0 + \tau_n$, which is the last component of $X_n$.
    Then, inductively define the sequence
    \begin{equation*}
        X_n := (G(Y_{\tau_n-}, \kappa_{\tau_n}, \eta_n), \kappa_{\tau_n}, R_{\tau_n-}, t_0+\tau_n), \quad n \ge 1.
    \end{equation*}
    By the piecewise formulation of $Y$, we can also write
    \begin{equation} \label{equation:definition-discretized-dynamics-2}
        X_n = (G(Y^{(n-1)}_{\tau_n - \tau_{n-1}}, \kappa_{\tau_n}, \eta_n), \kappa_{\tau_n}, R_{\tau_n-}, t_0+\tau_n), \quad n \ge 1.
    \end{equation}
    Since the fund's balance never exceeds $\bar{r}$ (see the proof of \Cref{lemma:functional-rep-of-R}), each $X_n$ takes values in the state space $\mathbb{X}$ of \Cref{definition:system-state-space}.
\end{definition}

\begin{definition}[Data History] \label{definition:data-history}
    Define the history of the above sequence up to index $n$ via
    \begin{equation} \label{equation:definition-data-history-1}
        H_n := (X_0, \xi_0, X_1, \xi_1, ..., X_{n-1}, \xi_{n-1}, X_n).
    \end{equation}
    Note the most current action $\xi_n$ is not included in $H_n$.
\end{definition}

\begin{definition}[Policies] \label{definition:policies}
    For $n \ge 0$, let $\mathbb{H}_n := (\mathbb{X} \times \mathbb{R}_+^d)^n \times \mathbb{X}$ be the space of histories $h_n = (x_0, a_0, \dots, x_{n-1}, a_{n-1}, x_n)$.
    A \emph{policy} is a sequence $\varpi = (\varpi_n)_{n \ge 0}$ of measurable functions $\varpi_n : \mathbb{H}_n \times [0, 1) \rightarrow \mathbb{R}_+^d$ such that $\varpi_n(h_n, \cdot)$ takes values in $A(x_n)$ for every $h_n \in \mathbb{H}_n$, where $x_n$ is the last component of $h_n$.
    Under the policy $\varpi$, the injection at time $\tau_n$ is
    \[
        \xi_n := \varpi_n(H_n, \chi_n), \qquad n \ge 0.
    \]
    We denote by $\Pi$ the set of all policies.
    A policy is \emph{deterministic} if no $\varpi_n$ depends on its last argument, and \emph{stationary deterministic} if, in addition, there is a measurable function $\varpi: \mathbb{X} \rightarrow \mathbb{R}_+^d$ with $\varpi(x) \in A(x)$ for all $x \in \mathbb{X}$ such that $\varpi_n(h_n, \cdot) \equiv \varpi(x_n)$ for all $n$ and $h_n$; we then identify the policy with the function $\varpi$ and write $\xi_n = \varpi(X_n)$.
\end{definition}

The sequence $X_0, \xi_0, X_1, \xi_1, \dots$ is constructed inductively: $X_0$ is given, and once $H_n$ is known, $\xi_n = \varpi_n(H_n, \chi_n)$ determines $Y^{(n)}_0$ and $R_{\tau_n}$, hence the dynamics on $[\tau_n, \tau_{n+1})$ and $X_{n+1}$. Throughout the paper, the injections are chosen by a policy $\varpi \in \Pi$.

\begin{remark} \label{remark:policies}
    Under a policy, the fund's manager observes the system only at signal times: when choosing the injection at $\tau_n$, the manager knows the past states $X_0, \dots, X_n$ and the past injections $\xi_0, \dots, \xi_{n-1}$, but not the path of $Y$ between signals.
    The variable $\chi_n$ allows the manager to randomize: since $\chi_n$ is independent of $H_n$, given $H_n = h_n$ the injection $\xi_n$ has distribution $\P\{\varpi_n(h_n, \chi_n) \in \cdot\}$, which is concentrated on $A(x_n)$.
    This sequence of stochastic kernels is a policy in the sense of \textcite{Feinberg}, which allows us to apply their results in \Cref{theorem:existence-of-optimal-policy}. We do not need the converse, namely that every policy of \textcite{Feinberg}, a sequence of stochastic kernels from histories to actions, can be written in the form of \Cref{definition:policies}: the stationary policy $\varpi^*$ of \Cref{theorem:existence-of-optimal-policy} attains the value function of \textcite{Feinberg} for the cost $\hat{c}$, so it is optimal also among all their policies (see the proof of that theorem).
\end{remark}

\begin{lemma} \label{lemma:functional-rep-of-Y}
    There exists a measurable function $\phi: \mathbb{R}^d \times \mathbb{R}_+ \times \mathscr{C}(\mathbb{R}_+, \mathbb{R}^d) \rightarrow \mathbb{R}^d$ such that
    \begin{equation} \label{equation:lemma-functional-rep-of-Y-1}
        Y^{(n)}_{\tau_{n+1}-\tau_n}= \phi(Y^{(n)}_0, \tau_{n+1}-\tau_n, B^{(n)}) \quad \text{a.s. } \forall n\in\mathbb{N}.
    \end{equation}
\end{lemma}
\begin{proof}
    Define $\phi(y, t, w) := \mathfrak{F}(y, w)(t)$, where $\mathfrak{F}$ is from \nCref{lemma:functional-solution-to-SDE}, which immediately implies \cref{equation:lemma-functional-rep-of-Y-1}.
    The measurability of $\phi$ follows from the measurability of $\mathfrak{F}$ and \Cref{lemma:measurability-of-evaluation-map}.
\end{proof}

\begin{lemma} \label{lemma:ODE-solution}
    There exists a measurable function $\Psi: \mathbb{R}_+ \times \mathbb{R}_+ \times \mathscr{C}(\mathbb{R}_+, \mathbb{R}^d) \rightarrow \mathbb{R}_+$ such that
    \begin{equation} \label{equation:lemma-ODE-solution-1}
        \mathfrak{R}_t := \Psi(r_0, t, y)
    \end{equation}
    satisfies
    \begin{equation} \label{equation:lemma-ODE-solution-2}
        \mathfrak{R}_t = r_0 + \int_0^t \Big(\rho \mathfrak{R}_s + \sum_{j=1}^d \delta_j (y^j_s \vee 0)\Big) \dd{s}.
    \end{equation}
\end{lemma}
\begin{proof}
    \Cref{equation:lemma-ODE-solution-2} is the classical ODE which can be solved explicitly via integrating factor technique.
    Its solution is
    \begin{equation} \label{equation:lemma-ODE-solution-3}
        \mathfrak{R}_t = e^{\rho t} r_0 + \sum_{j=1}^d \delta_j \int_0^t e^{\rho (t-s)}(y^j_s \vee 0) \dd{s}.
    \end{equation}
    Define auxiliary function $I_{j,m}:\mathbb{R}_+ \times \mathscr{C}(\mathbb{R}_+, \mathbb{R}^d) \rightarrow \mathbb{R}_+$ via
    \begin{equation} \label{equation:lemma-ODE-solution-4}
        I_{j,m}(t, y) := \sum_{i=0}^{m-1} \frac{t}{m}\cdot e^{\rho (t-\frac{ti}{m})}\left[\mathfrak{e}\left(\frac{ti}{m}, \pi_j(y) \right) \vee 0 \right],
    \end{equation}
    where evaluation map $\mathfrak{e}$ is defined in \Cref{lemma:measurability-of-evaluation-map} and $\pi_j$ projects $y$ to its j\textsuperscript{th} coordinate.
    Because $y$ is continuous, the above Riemann sum converges to the corresponding integral.
    Hence,
    \begin{equation} \label{equation:lemma-ODE-solution-5}
        I_{j,m}(t, y) \overset{m \rightarrow \infty}{\longrightarrow} \int_0^t e^{\rho (t-s)}(y^j_s \vee 0) \dd{s} =: I_j(t, y) \quad \text{pointwise}.
    \end{equation}
    The joint measurability of $I_j$ follows from the measurability of $\mathfrak{e}$ and $\pi_j$, and the fact that convergence preserves measurability.
    Finally, by \cref{equation:lemma-ODE-solution-3}, we have
    \begin{equation} \label{equation:lemma-ODE-solution-6}
        \mathfrak{R}_t = e^{\rho t}r_0 + \sum_{j=1}^d \delta_j I_j(t, y) =: \Psi(r_0, t, y),
    \end{equation}
    where $\Psi$'s measurability follows from $I_j$'s.
\end{proof}

\begin{lemma} \label{lemma:functional-rep-of-R}
    There exists a measurable function $\psi: \mathbb{R}^d \times \mathbb{R}_+ \times \mathbb{R}_+ \times \mathscr{C}(\mathbb{R}_+, \mathbb{R}^d) \rightarrow \mathbb{R}_+$ such that
    \begin{equation} \label{equation:lemma-functional-rep-of-R-1}
        R_{\tau_{n+1}-} = \psi(Y^{(n)}_0, R_{\tau_n}, \tau_{n+1}-\tau_n, B^{(n)}) \quad \text{a.s. } \forall n \in \mathbb{N}.
    \end{equation}
\end{lemma}
\begin{proof}
    The piecewise representation of process $R$ gives
    \begin{equation} \label{equation:lemma-functional-rep-of-R-2}
        R_{\tau_{n+1}-} = R^{(n)}_{\tau_{n+1}-\tau_n},
    \end{equation}
    where $R^{(n)}$ follows
    \begin{equation} \label{equation:lemma-functional-rep-of-R-3}
        R^{(n)}_t = \Psi(R_{\tau_n}, t, Y^{(n)}) \wedge \bar{r}
    \end{equation}
    with $\Psi$ defined in the previous lemma.
    This holds because $0 \le R_{\tau_n} \le \bar{r}$ for all $n$, by induction: $R_{\tau_0} = r_0 - ||\xi_0||_1 \in [0, \bar{r}]$ because $\xi_0 \in A(X_0)$; and if $0 \le R_{\tau_n} \le \bar{r}$, then $R$ follows \cref{equation:money-reserve-1} on $[\tau_n, \tau_{n+1})$ starting from $R_{\tau_n}$, so $R_{\tau_{n+1}-} = \Psi(R_{\tau_n}, \tau_{n+1}-\tau_n, Y^{(n)}) \wedge \bar{r} \in [0, \bar{r}]$ (see the discussion after \cref{equation:money-reserve-1}), and $R_{\tau_{n+1}} = R_{\tau_{n+1}-} - ||\xi_{n+1}||_1 \in [0, \bar{r}]$ because $\xi_{n+1} \in A(X_{n+1})$.
    Now, since $Y^{(n)} = \mathfrak{F}(Y^{(n)}_0, B^{(n)})$ by \Cref{lemma:functional-solution-to-SDE}, the existence of a measurable function $\psi$ follows from the measurability of $\Psi$ and $\mathfrak{F}$.
\end{proof}

\begin{theorem}[Markov Property] \label{theorem:markov-property}
    Under any policy $\varpi \in \Pi$, for every $n \ge 0$ and any Borel set $\Gamma \subset \mathbb{X}$,
    \begin{equation} \label{equation:theorem-markov-property-1}
        \P\{X_{n+1} \in \Gamma | H_n, \xi_n\} = \P\{X_{n+1} \in \Gamma | X_n, \xi_n\}.
    \end{equation}
\end{theorem}
\begin{proof}
    First, recall $X_{n+1}$ is given by
    \begin{equation*} \label{equation:theorem-markov-property-2}
        \begin{pmatrix}
            G(Y^{(n)}_{\tau_{n+1}-\tau_n}, \kappa_{\tau_{n+1}}, \eta_{n+1}) \\
            \kappa_{\tau_{n+1}} \\
            R_{\tau_{n+1}-} \\
            t_0 + \tau_{n+1}
        \end{pmatrix} = \begin{pmatrix}
            G(\phi(Y^{(n)}_0, \tau_{n+1}-\tau_n, B^{(n)}), \kappa_{\tau_{n+1}}, \eta_{n+1}) \\
            \kappa_{\tau_{n+1}} \\
            \psi(Y^{(n)}_0, R_{\tau_n}, \tau_{n+1}-\tau_n, B^{(n)}) \\
            t_0 + \tau_n + (\tau_{n+1}-\tau_n)
        \end{pmatrix},
    \end{equation*}
    where we used \Cref{lemma:functional-rep-of-Y} and \Cref{lemma:functional-rep-of-R}.
    We could further substitute in $Y^{(n)}_0 = X_{n, 1}+\xi_n$ and $R_{\tau_n} = X_{n, 3}-\sum_{j=1}^d \xi_{n,j}$, where $X_{n, l}$ means the $l^\text{th}$ component of $X_n$, and $\xi_{n,j}$ means the $j\text{th}$ coordinate of $\xi_n$.
    Define our first auxiliary function $g_1: \text{Gr}(A) \times \mathbb{R}_+ \times \mathscr{K} \times [0, 1)^{|\mathscr{K}| \times d} \times \Omega \rightarrow [0, 1]$ by mapping $(x, a, s, k', u, \omega)$, with notation $x=(y, k, r, t) \in \mathbb{X}$, to
    \begin{equation*} \label{equation:theorem-markov-property-3}
        1_\Gamma\left(G\big(\phi(y+a, s, B^{(n)}(\omega)), k', u\big), k', \psi\bigg(y+a, r-\sum_{j=1}^d a_j, s, B^{(n)}(\omega)\bigg), t+s \right).
    \end{equation*}
    Function $g_1$ is jointly measurable (use the usual discrete $\sigma$-field for $\mathscr{K}$) by the measurability of $\phi, \psi, G$. On $\text{Gr}(A)$ we have $r - \sum_{j=1}^d a_j \ge 0$, so the second argument of $\psi$ lies in its domain $\mathbb{R}_+$; this is why $g_1$ and the functions below are defined on $\text{Gr}(A)$. Since $\xi_n \in A(X_n)$, we have $(X_n, \xi_n) \in \text{Gr}(A)$.
    Define the second auxiliary function $g_2: \text{Gr}(A) \times \mathbb{R}_+ \times \mathscr{K} \times [0, 1)^{|\mathscr{K}|\times d} \rightarrow [0, 1]$ via
    \begin{equation*} \label{equation:theorem-markov-property-4}
        g_2(x, a, s, k', u) := \mathbb{E}g_1(x, a, s, k', u, \cdot).
    \end{equation*}
    We make an important observation that is essential for the next theorem.
    \begin{quote}
        (Observation 1) When computing $g_2$, because $B^{(n)} = B$ in distribution, if we had used $B$ in place of $B^{(n)}$ in $g_1$'s definition, we would still get the same function $g_2$.
    \end{quote}
    Since $\mathcal{G} \subseteq \mathcal{F}_0$, the random variables $H_n$, $\xi_n$, $\tau_{n+1}-\tau_n$, $\kappa_{\tau_{n+1}}$ and $\eta_{n+1}$ are $\mathcal{F}_{\tau_n}$-measurable, and $B^{(n)}$ is independent of $\mathcal{F}_{\tau_n}$.
    By \nCref{lemma:freezing-lemma}, applied with $\mathcal{H} := \sigma(B^{(n)})$ and $\mathcal{D} := \sigma(H_n, \xi_n, \tau_{n+1}-\tau_n, \kappa_{\tau_{n+1}}, \eta_{n+1}) \subseteq \mathcal{F}_{\tau_n}$,
    \begin{equation*} \label{equation:theorem-markov-property-5}
        \P\{X_{n+1} \in \Gamma | H_n, \xi_n, \tau_{n+1}-\tau_n, \kappa_{\tau_{n+1}}, \eta_{n+1}\} = g_2(X_n, \xi_n, \tau_{n+1}-\tau_n, \kappa_{\tau_{n+1}}, \eta_{n+1}).
    \end{equation*}

    We now compute the expectation of $g_2$ conditioning on $(H_n, \xi_n, \tau_{n+1}-\tau_n, \kappa_{\tau_{n+1}})$.
    Because $\eta_{n+1}$ is independent of $(H_n, \xi_n, \tau_{n+1}-\tau_n, \kappa_{\tau_{n+1}})$, which are functions of $B$, $(E_i, U_i)_{i \le n+1}$, $(\eta_i)_{i \le n}$ and $(\chi_i)_{i \le n}$ (\Cref{assumption:standing}), \Cref{lemma:freezing-lemma} gives
    \begin{equation*} \label{equation:theorem-markov-property-6}
        \begin{split}
            &\quad\,\,\mathbb{E}\{g_2(X_n, \xi_n, \tau_{n+1}-\tau_n, \kappa_{\tau_{n+1}}, \eta_{n+1}) | H_n, \xi_n, \tau_{n+1}-\tau_n, \kappa_{\tau_{n+1}}\} \nonumber\\
            &= \int_{[0, 1)^{|\mathscr{K}|\times d}} g_2(X_n, \xi_n, \tau_{n+1}-\tau_n, \kappa_{\tau_{n+1}}, u) \dd{\mu_\eta(u)},
        \end{split}
    \end{equation*}
    where $\mu_\eta$ is the i.i.d. distribution of $\eta_n$.
    Then there exists a measurable function $g_3: \text{Gr}(A) \times \mathbb{R}_+ \times \mathscr{K} \rightarrow [0, 1]$ (our third auxiliary function) such that the above quantity equals to $g_3(X_n, \xi_n, \tau_{n+1}-\tau_n, \kappa_{\tau_{n+1}})$ a.s.

    The last step is to compute the expectation of $g_3$ conditioning on $(H_n, \xi_n)$.
    We invoke \Cref{lemma:more-on-signal-interval-distribution} and \Cref{remark:further-on-signal-interval-distribution} to obtain
    \begin{equation*} \label{equation:theorem-markov-property-7}
        \begin{split}
            &\quad\,\,\mathbb{E}\{g_3(X_n, \xi_n, \tau_{n+1}-\tau_n, \kappa_{\tau_{n+1}}) | H_n, \xi_n\} \\
            &= \sum_{k' \in \mathscr{K}} \int_0^\infty g_3(X_n, \xi_n, s, k') \lambda_{\kappa_{\tau_n}, k'}e^{-s\lambda_{\kappa_{\tau_n}}} \dd{s},
        \end{split}
    \end{equation*}
    and we name this new quantity $g_4(X_n, \xi_n)$ for our fourth auxiliary measurable function.
    We have now (implicitly) computed the conditional probability $\P\{X_{n+1} \in \Gamma | H_n, \xi_n\} = g_4(X_n, \xi_n)$. Since $g_4(X_n, \xi_n)$ is $\sigma(X_n, \xi_n)$-measurable, the tower property gives $\P\{X_{n+1} \in \Gamma | X_n, \xi_n\} = g_4(X_n, \xi_n)$ as well, which is \cref{equation:theorem-markov-property-1}.
\end{proof}

\begin{theorem} \label{theorem:transition-kernel}
    There exists a stochastic kernel $Q$ on $\mathbb{X}$ given $\text{Gr}(A)$ such that
    \begin{equation} \label{equation:theorem-transition-kernel-1}
        \P\{X_{n+1} \in \Gamma | X_n, \xi_n\} = Q(\Gamma | X_n, \xi_n)
    \end{equation}
    for all $n$, all policies $\varpi \in \Pi$ and Borel set $\Gamma$.
\end{theorem}
\begin{proof}
    In the proof of \Cref{theorem:markov-property}, the functions $g_2$, $g_3$ and $g_4$ do not depend on $n$. For $g_2$ this is (Observation 1). For $g_3$, all the $\eta_n$ have the same distribution $\mu_\eta$. For $g_4$, the conditional distribution of $(\tau_{n+1}-\tau_n, \kappa_{\tau_{n+1}})$ given $(H_n, \xi_n)$ is given by the same formula for every $n$ (\Cref{remark:more-on-signal-interval-distribution,remark:further-on-signal-interval-distribution}). Nor do these functions depend on the policy $\varpi$, which enters only through $\xi_n$. Hence we have
    \begin{equation} \label{equation:theorem-transition-kernel-2}
        \P\{X_{n+1} \in \Gamma | X_n, \xi_n\} = g_4^\Gamma(X_n, \xi_n) \quad \forall n,
    \end{equation}
    where $g_4^\Gamma$ is the function $g_4$ defined in the previous theorem's proof; we add the superscript here to emphasize it is induced from the given set $\Gamma$.
    Now, for any Borel set $\Gamma \subset \mathbb{X}$, $x=(y, k, r, t) \in \mathbb{X}$ and $a\in A(x)$, define
    \begin{equation} \label{equation:theorem-transition-kernel-3}
        Q(\Gamma | x, a) := g_4^\Gamma(x, a).
    \end{equation}
    Claim this is the desired stochastic kernel.
    Firstly, \cref{equation:theorem-transition-kernel-1} is immediate by eq. (\ref{equation:theorem-transition-kernel-2})(\ref{equation:theorem-transition-kernel-3}).
    The task now is to check that $Q$ is a stochastic kernel.

    Fix an arbitrary $x=(y, k, r, t) \in \mathbb{X}$ and $a \in A(x)$.
    Define the random vector
    \begin{equation} \label{equation:theorem-transition-kernel-4}
        \begin{split}
            &Z^{(x, a)} := \left(G(Y^{y+a}_S, K, \eta), K, \vphantom{\int_0^S} \right.\\
            &\hspace{5em} \left. \bar{r} \wedge \bigg[e^{\rho S}(r-||a||_1) + \sum_{j=1}^d \delta_j \int_0^S e^{\rho (S-u)}(Y^{y+a, j}_u \vee 0) \dd{u} \bigg], t+S \right)
        \end{split}
    \end{equation}
    with $(S, K)$ being the time until the next market signal and the mark of that signal, given that the most recent signal has mark $k$ (that is, $S$ is exponential with rate $\lambda_k$, $K$ is independent of $S$, and $\P\{K = k'\} = \lambda_{k,k'}/\lambda_k$ for $k' \in \mathscr{K}$), $Y^{y+a}$ being the solution of \cref{equation:wealth-dynamics-3} with initial condition $y+a$, and $\eta$ following the same distribution as $\eta_n$'s. For instance, we may take $S := E_1/\lambda_k$, $K := F(k, U_1)$ and $\eta := \eta_1$, and let $Y^{y+a}$ be driven by $B$; then $S$, $K$, $\eta$ and $B$ are independent (\Cref{assumption:standing}).
    By going through the computations in the previous theorem again, we see that
    \begin{equation} \label{equation:theorem-transition-kernel-5}
        Q(\cdot | x, a) = \P\{Z^{(x, a)} \in \cdot\},
    \end{equation}
    which confirms that $Q(\cdot | x, a)$ is a probability measure for any given $(x, a)$.
    Next, fix an arbitrary Borel set $\Gamma$.
    By \cref{equation:theorem-transition-kernel-3} and $g_4^\Gamma$ being measurable (use the usual subspace Borel $\sigma$-algebra), $Q(\Gamma | \cdot)$ is a measurable function as required for a stochastic kernel.
\end{proof}

\begin{theorem}[Weak Continuity] \label{theorem:weak-continuity}
    Take any sequence $(x_m, a_m) \rightarrow (x, a)$ in $\text{Gr}(A)$ (\Cref{definition:feasible-action-set}), then transition kernel $Q(\cdot | x_m, a_m)$ converges to $Q(\cdot | x, a)$ weakly.
\end{theorem}
\begin{proof}
    By the previous theorem, $Q(\cdot | x, a)$ is the probability measure of $Z := Z^{(x, a)}$, and $Q(\cdot | x_m, a_m)$ is the probability measure of $Z_m := Z^{(x_m, a_m)}$.
    In what follows, we use notation $x = (y, k, r, t)$ and $x_m = (y_m, k_m, r_m, t_m)$.
    Because $\mathscr{K}$ is a discrete space, $x_m \rightarrow x$ implies $k_m = k$ for all large enough $m$.
    This means we can take $Z, Z_m$ to share the same random variables $S, K$, $\eta$ and the same Brownian motion $B$ (refer to \cref{equation:theorem-transition-kernel-4} for the meaning of $S, K$ and $\eta$) without loss of generality: with the choice made after \cref{equation:theorem-transition-kernel-4}, $Z_m$ and $Z$ are built from the same $S$, $K$, $\eta$ and $B$ once $k_m = k$.
    Now, we show $Z_m$ converges to $Z$ in probability, a stronger statement than the needed convergence in distribution.
    We use the fact that random vectors converge in probability if and only if the convergence holds component-wise.
    The second components of $Z_m$ and $Z$ coincide (both equal $K$), and the fourth components $t_m + S$ converge to $t + S$ a.s., since $t_m \rightarrow t$.

    We next focus on the first component.
    Observe from the definition of $G$ that
    \begin{equation*} \label{equation:theorem-weak-continuity-2}
        ||G(Y^{y_m+a_m}_S, K, \eta) - G(Y^{y+a}_S, K, \eta)|| \le ||Y^{y_m+a_m}_S - Y^{y+a}_S||.
    \end{equation*}
    Invoke \Cref{lemma:SDE-convergence}, for any $\epsilon' > 0$, we have
    \begin{equation*} \label{equation:theorem-weak-continuity-3}
        \begin{split}
            &\quad\,\,\P\{||Y^{y_m+a_m}_S - Y^{y+a}_S|| > \epsilon'\} \\
            &\le \P\left\{\sup_{0\le u \le s} ||Y^{y_m+a_m}_u - Y^{y+a}_u|| > \epsilon', S \le s\right\} + \P\{S > s\} \\
            &\le \frac{C||y_m+a_m-(y+a)||^2}{\epsilon'^2} + \P\{S > s\},
        \end{split}
    \end{equation*}
    where $C$ depends only on $s$ and the Lipschitz constants of the coefficient functions $b, \sigma$ from $Y$'s SDE.
    Now, for any $\epsilon > 0$, because $S < \infty$ a.s., we could pick $s$ large enough such that $\P\{S > s\} < \epsilon$.
    With this $s$ now fixed, constant $C$ is also settled.
    Then, pass $(y_m, a_m) \rightarrow (y, a)$, the right hand side above can be made less than $2\epsilon$.

    Now, we focus on the third component of $Z_m$.
    Since $z \mapsto \bar{r} \wedge z$ is 1-Lipschitz, it suffices to show that the expression in square brackets in \cref{equation:theorem-transition-kernel-4} converges in probability.
    The term $e^{\rho S}(r_m-||a_m||_1)$ converges to $e^{\rho S}(r-||a||_1)$ a.s. since $r_m \rightarrow r$ and $a_m \rightarrow a$.
    For the integral term, we use the same trick of breaking into the cases $S > s$ and $S \le s$. Since $|z \vee 0 - z' \vee 0| \le |z - z'|$, on the event $\{S \le s\}$ we have
    \[
    \begin{split}
        \bigg|\sum_{j=1}^d \delta_j \int_0^S e^{\rho (S-u)}\big[(Y^{y_m+a_m, j}_u \vee 0) - (Y^{y+a, j}_u \vee 0)\big] \dd{u}\bigg|  \\
        \le s e^{\rho s} \sum_{j=1}^d \delta_j \sup_{0 \le u \le s} ||Y^{y_m+a_m}_u - Y^{y+a}_u||,
    \end{split}
    \]
    so the convergence in probability follows from \Cref{lemma:SDE-convergence}, as for the first component.

    Hence, we have established that $Z_m \rightarrow Z$ in probability, which gives the desired weak convergence of $Q$.
\end{proof}

We now define the objective function. Fix a policy $\varpi \in \Pi$ and an initial state $x_0 \in \mathbb{X}$, let $X_0 := x_0$, and let the injections $\xi_n$ be chosen by $\varpi$ (\Cref{definition:policies}). The expected total cost is
\begin{equation} \label{equation:objective-function-1}
    V(\varpi, x_0) := \mathbb{E}\sum_{n=0}^\infty c(X_n, \xi_n),
\end{equation}
and its associated $N$-stage total cost is
\begin{equation} \label{equation:objective-function-2}
    V_N(\varpi, x_0) := \mathbb{E}\sum_{n=0}^{N-1} c(X_n, \xi_n);
\end{equation}
their optimal value functions are given by
\begin{equation} \label{equation:objective-function-3}
    V^*(x_0) := \inf_{\varpi \in \Pi} V(\varpi, x_0), \text{ and } V_N^*(x_0) := \inf_{\varpi \in \Pi} V_N(\varpi, x_0).
\end{equation}
A technical subtlety occurs: we need to know if the series above converges, let alone whether it produces a finite expectation.
This is an issue because cost function $c$ takes both positive and negative values.
A natural step is to find the set of policies that make $V(\varpi, x_0)$ properly defined.
Fortunately, the following few estimating results tell us that, when the discount rate $\gamma$ is large enough (\Cref{theorem:one-stage-cost-estimate}), \textit{all} policies are allowed.
Furthermore, they imply an interesting convergence rate result between $V^*_N$ and $V^*$.

\begin{lemma} \label{lemma:cost-bound}
    Let $x=(y, k, r, t) \in \mathbb{X}$ be the system's state value and $a \in \mathbb{R}_+^d$ be any feasible action value, i.e. $||a||_1 \le r$.
    Then there exists a constant $C_\text{cost}$ independent of $x, a$ satisfying
    \begin{equation} \label{equation:lemma-cost-bound-1}
        |c(x, a)| \le C_\text{cost}e^{-\gamma t}\left(\sqrt{1 + ||y||^2} + r \right).
    \end{equation}
\end{lemma}
\begin{proof}
    Since $a_j \ge 0$ and $l_j \ge 0$ (\Cref{assumption:standing}), we have $(l_j - y_j - a_j)^+ \le (l_j - y_j)^+ \le l_j + |y_j|$, so the linear growth of $f$ gives $0 \le f(l_j - y_j - a_j) \le C_f(1 + l_j + |y_j|)$.
    Moreover, $0 \le ||a||_1 \le r$, so $1 \le 1 + r - ||a||_1 \le 1 + r$ and hence $0 \le (1 + r - ||a||_1)^\nu \le 1 + r$, because $\nu \in (0, 1)$.
    By the triangle inequality and these bounds,
    \begin{equation*}
        \begin{split}
            e^{\gamma t}|c(x, a)| &\le \sum_{j=1}^d w_j f(l_j - y_j - a_j) + w_A||a||_1 + w_R(1 + r - ||a||_1)^\nu \\
            &\le C_f\sum_{j=1}^d w_j(1 + l_j) + C_f\sum_{j=1}^d w_j|y_j| + w_A r + w_R(1 + r) \\
            &\le \Big[C_f\sum_{j=1}^d w_j(1 + l_j) + w_R\Big] + C_f\Big(\sum_{j=1}^d w_j^2\Big)^{1/2}||y|| + (w_A + w_R)\, r,
        \end{split}
    \end{equation*}
    where the last step uses the Cauchy--Schwarz inequality.
    Since $1 \le \sqrt{1 + ||y||^2}$ and $||y|| \le \sqrt{1 + ||y||^2}$, \cref{equation:lemma-cost-bound-1} holds with
    \[
        C_\text{cost} := C_f\Big[\sum_{j=1}^d w_j(1 + l_j) + \Big(\sum_{j=1}^d w_j^2\Big)^{1/2}\Big] + w_A + 2w_R.
    \]
\end{proof}
\begin{remark}
    Notice the bound above is independent of feasible actions $a$. This is essential when we bound the objective function $V$ under arbitrary policies.
\end{remark}

\begin{theorem}[One Stage Cost Estimate] \label{theorem:one-stage-cost-estimate}
    Let $\Lambda_\text{max}(\cdot)$ denote the maximum eigenvalue of a given matrix.
    Let
    \[
        \zeta := \max\Big(\Lambda_\text{max}\big(\tfrac{\Theta+\Theta^\intercal}{2}\big) + \tfrac{1}{2}\max_{j=1,\dots,d} \sigma_j^2,\ 0\Big), \qquad q := \frac{\bar{\lambda}}{\bar{\lambda}+\gamma-\zeta},
    \]
    where the matrix $\Theta$ is defined in \cref{equation:wealth-dynamics-2} and $\bar{\lambda}$ in \cref{equation:signal-total-rate}.
    If the discount rate $\gamma$ satisfies $\gamma > \zeta$, then $q \in (0, 1)$ and, for every policy $\varpi \in \Pi$, every initial state $x_0 = (y_0, k_0, r_0, t_0) \in \mathbb{X}$ and every $n \ge 0$,
    \begin{equation} \label{equation:theorem-one-stage-cost-estimate-1}
        \mathbb{E}|c(X_n, \xi_n)| \le C_\text{cost}\, e^{-\gamma t_0} \Big(\sqrt{1 + ||y_0||^2} + (n+1)\bar{r}\Big) q^n,
    \end{equation}
    where $C_\text{cost}$ is the constant of \Cref{lemma:cost-bound}. In particular, $\mathbb{E}|c(X_n, \xi_n)| = \mathscr{O}(nq^n)$ uniformly in $\varpi \in \Pi$, where $\mathscr{O}$ is the usual big-O notation.
\end{theorem}
\begin{proof}
    Fix $\varpi \in \Pi$ and $x_0 \in \mathbb{X}$. Throughout this proof, write $\tilde{Y}_n := X_{n,1}$ for the first component of $X_n$ (the members' wealth right after the $n$-th signal, before the injection) and $S_n := \tau_{n+1} - \tau_n$, so that $\tilde{Y}_0 = y_0$, $Y^{(n)}_0 = \tilde{Y}_n + \xi_n$ and $\tilde{Y}_{n+1} = G(Y^{(n)}_{S_n}, \kappa_{\tau_{n+1}}, \eta_{n+1})$.
    Let also $p := \bar{\lambda}/(\bar{\lambda} + \gamma)$; since $\zeta \ge 0$, we have $0 < p \le q < 1$.

    \emph{Step 1: reduction.} Since $X_{n,3} = R_{\tau_n-} \le \bar{r}$ and $X_{n,4} = t_0 + \tau_n$, \Cref{lemma:cost-bound} gives
    \begin{equation} \label{equation:one-stage-proof-1}
        \mathbb{E}|c(X_n, \xi_n)| \le M_n := C_\text{cost}\,\mathbb{E}\Big\{e^{-\gamma(t_0+\tau_n)}\Big(\sqrt{1 + ||\tilde{Y}_n||^2} + \bar{r}\Big)\Big\}.
    \end{equation}

    \emph{Step 2: the wealth term.} Each component of $G(y, k, u)$ satisfies $|G(y,k,u)_j| \le |y_j|$, so $||\tilde{Y}_{n+1}|| \le ||Y^{(n)}_{S_n}||$, and $Y^{(n)}_{S_n} = \phi(Y^{(n)}_0, S_n, B^{(n)})$ a.s.\ by \Cref{lemma:functional-rep-of-Y}.
    The random variables $\tau_n$, $S_n$ and $Y^{(n)}_0$ are $\mathcal{F}_{\tau_n}$-measurable ($\tau_n$ and $S_n$ are even $\mathcal{F}_0$-measurable), and $B^{(n)}$ is independent of $\mathcal{F}_{\tau_n}$.
    Hence \Cref{lemma:freezing-lemma}, applied to a nonnegative function, gives
    \[
        \mathbb{E}\Big\{e^{-\gamma S_n}\sqrt{1 + ||Y^{(n)}_{S_n}||^2}\,\Big|\,\mathcal{F}_{\tau_n}\Big\} = h(Y^{(n)}_0, S_n), \quad h(y, s) := e^{-\gamma s}\,\mathbb{E}\sqrt{1 + ||Y^y_s||^2},
    \]
    where $Y^y$ is the solution of \cref{equation:wealth-dynamics-3} with $Y^y_0 = y$, since $B^{(n)}$ has the same law as $B$.
    For $b(y) = \Theta y$ and $\sigma(y) = \text{diag}(\sigma_1 y_1, \dots, \sigma_d y_d)$,
    \[
        y^\intercal b(y) + \tfrac{1}{2}||\sigma(y)||^2 = y^\intercal \tfrac{\Theta + \Theta^\intercal}{2}\, y + \tfrac{1}{2}\sum_{j=1}^d \sigma_j^2 y_j^2 \le \zeta||y||^2 \le \zeta(1 + ||y||^2),
    \]
    so \Cref{lemma:SDE-estimate} and Jensen's inequality give $\mathbb{E}\sqrt{1 + ||Y^y_s||^2} \le e^{\zeta s}\sqrt{1 + ||y||^2}$, that is, $h(y, s) \le e^{-(\gamma - \zeta)s}\sqrt{1 + ||y||^2}$.
    Finally, $y \mapsto \sqrt{1 + ||y||^2}$ is 1-Lipschitz for the Euclidean norm, and $||\xi_n|| \le ||\xi_n||_1 \le X_{n,3} \le \bar{r}$ because $\xi_n \in A(X_n)$, so $\sqrt{1 + ||Y^{(n)}_0||^2} \le \sqrt{1 + ||\tilde{Y}_n||^2} + \bar{r}$.
    Since $e^{-\gamma(t_0+\tau_n)}$ is $\mathcal{F}_{\tau_n}$-measurable, combining these bounds yields
    \begin{equation} \label{equation:one-stage-proof-2}
        \mathbb{E}\Big\{e^{-\gamma(t_0+\tau_{n+1})}\sqrt{1 + ||\tilde{Y}_{n+1}||^2}\Big\} \le \mathbb{E}\Big\{e^{-\gamma(t_0+\tau_n)}e^{-(\gamma-\zeta)S_n}\Big(\sqrt{1 + ||\tilde{Y}_n||^2} + \bar{r}\Big)\Big\}.
    \end{equation}

    \emph{Step 3: the waiting time.} By \Cref{remark:more-on-signal-interval-distribution,remark:further-on-signal-interval-distribution} (summed over $k'$), given $(H_n, \xi_n)$ the waiting time $S_n$ is exponential with rate $\lambda_{\kappa_{\tau_n}}$, so for every $\vartheta > 0$,
    \[
        \mathbb{E}\{e^{-\vartheta S_n} | H_n, \xi_n\} = \frac{\lambda_{\kappa_{\tau_n}}}{\lambda_{\kappa_{\tau_n}} + \vartheta} \le \frac{\bar{\lambda}}{\bar{\lambda} + \vartheta},
    \]
    because $\lambda \mapsto \lambda/(\lambda + \vartheta)$ is increasing and $\lambda_{\kappa_{\tau_n}} \le \bar{\lambda}$.
    Since $\tau_n$ and $\tilde{Y}_n$ are functions of $X_n$, hence of $H_n$, conditioning the right-hand side of \cref{equation:one-stage-proof-2} on $(H_n, \xi_n)$ with $\vartheta = \gamma - \zeta > 0$ gives
    \[
        \mathbb{E}\Big\{e^{-\gamma(t_0+\tau_{n+1})}\sqrt{1 + ||\tilde{Y}_{n+1}||^2}\Big\} \le q\,\mathbb{E}\Big\{e^{-\gamma(t_0+\tau_n)}\Big(\sqrt{1 + ||\tilde{Y}_n||^2} + \bar{r}\Big)\Big\} = \frac{q M_n}{C_\text{cost}}.
    \]
    In the same way, with $\vartheta = \gamma$, $\mathbb{E}e^{-\gamma\tau_{n+1}} \le p\,\mathbb{E}e^{-\gamma\tau_n}$, so $\mathbb{E}e^{-\gamma\tau_n} \le p^n \le q^n$ for all $n$, since $\tau_0 = 0$.

    \emph{Step 4: recursion.} By Steps 2 and 3,
    \[
        M_{n+1} \le q M_n + C_\text{cost}\,\bar{r}\, e^{-\gamma t_0}\,\mathbb{E}e^{-\gamma\tau_{n+1}} \le q M_n + C_\text{cost}\,\bar{r}\, e^{-\gamma t_0} q^{n+1}.
    \]
    Since $M_0 = C_\text{cost} e^{-\gamma t_0}(\sqrt{1 + ||y_0||^2} + \bar{r})$, induction on $n$ gives $M_n \le C_\text{cost} e^{-\gamma t_0}(\sqrt{1 + ||y_0||^2} + (n+1)\bar{r}) q^n$, and \cref{equation:theorem-one-stage-cost-estimate-1} follows from \cref{equation:one-stage-proof-1}.
\end{proof}
\begin{remark}
    Like the bound in \Cref{lemma:cost-bound}, the bound in \cref{equation:theorem-one-stage-cost-estimate-1} does not depend on the injections, and hence it holds uniformly over all policies. When $\zeta > 0$, the rate can be improved to $\mathscr{O}(q^n)$: with $p := \bar{\lambda}/(\bar{\lambda}+\gamma) < q$, keeping track separately of the terms that come from the fund in the proof above gives $\mathbb{E}|c(X_n, \xi_n)| \le C_\text{cost}\, e^{-\gamma t_0}\big(\sqrt{1 + ||y_0||^2} + \bar{r} + \bar{r}\, q/(q - p)\big) q^n$. The exponential rate $q$ is the same, and the constant blows up as $\zeta \rightarrow 0$ (when $\zeta = 0$, $p = q$), so we work with \cref{equation:theorem-one-stage-cost-estimate-1}, which covers both cases.
\end{remark}

We now discuss two implications of the above estimates.

\begin{corollary} \label{corollary:V-defined}
    Under the assumption of \Cref{theorem:one-stage-cost-estimate}, for every policy $\varpi \in \Pi$ and every $x_0 \in \mathbb{X}$ we have $\mathbb{E}\sum_{n=0}^\infty |c(X_n, \xi_n)| < \infty$. In particular, the series $\sum_n c(X_n, \xi_n)$ converges absolutely a.s., $V(\varpi, x_0)$ and $V_N(\varpi, x_0)$ are finite for every $N$, and $V(\varpi, x_0) = \sum_{n=0}^\infty \mathbb{E}c(X_n, \xi_n)$.
\end{corollary}
\begin{proof}
    By monotone convergence and \cref{equation:theorem-one-stage-cost-estimate-1},
    \[
        \mathbb{E}\sum_{n=0}^\infty |c(X_n, \xi_n)| = \sum_{n=0}^\infty \mathbb{E}|c(X_n, \xi_n)| \le C_\text{cost}\, e^{-\gamma t_0}\sum_{n=0}^\infty \Big(\sqrt{1 + ||y_0||^2} + (n+1)\bar{r}\Big) q^n < \infty,
    \]
    since $q < 1$. Hence $\sum_n |c(X_n, \xi_n)| < \infty$ a.s., and the other statements follow from dominated convergence, with $\sum_n |c(X_n, \xi_n)|$ as the dominating function.
\end{proof}

This resolves the question of whether $V(\varpi, x_0)$ is well-defined.
Next, we derive the convergence rate of $V_N^*$ to $V^*$.

\begin{corollary} \label{corollary:V_N-convergence-rate}
    Under the assumption of \Cref{theorem:one-stage-cost-estimate}, for every $x_0 = (y_0, k_0, r_0, t_0) \in \mathbb{X}$ and every $N \ge 1$,
    \begin{equation} \label{equation:corollary-V_N-convergence-rate-1}
        \sup_{\varpi \in \Pi} |V_N(\varpi, x_0) - V(\varpi, x_0)| \le \Delta_N(x_0) := C_\text{cost}\, e^{-\gamma t_0} q^N \bigg[\frac{\sqrt{1 + ||y_0||^2} + (N+1)\bar{r}}{1-q} + \frac{q\,\bar{r}}{(1-q)^2}\bigg],
    \end{equation}
    with $q$ as in \Cref{theorem:one-stage-cost-estimate}. Consequently, $V^*(x_0)$ and $V_N^*(x_0)$ are finite, $|V_N^*(x_0) - V^*(x_0)| \le \Delta_N(x_0)$, and in particular $|V_N^*(x_0) - V^*(x_0)| = \mathscr{O}(Nq^N)$.
\end{corollary}
\begin{proof}
    Fix $\varpi \in \Pi$. By \Cref{corollary:V-defined}, $V(\varpi, x_0) - V_N(\varpi, x_0) = \sum_{n \ge N} \mathbb{E}c(X_n, \xi_n)$, so by \Cref{theorem:one-stage-cost-estimate},
    \[
        \begin{split}
            |V(\varpi, x_0) & - V_N(\varpi, x_0)|  \le \sum_{n \ge N} \mathbb{E}|c(X_n, \xi_n)| \\
            & \le C_\text{cost}\, e^{-\gamma t_0}\sum_{m=0}^\infty \Big(\sqrt{1 + ||y_0||^2} + (N+1+m)\bar{r}\Big) q^{N+m} = \Delta_N(x_0),
        \end{split}
    \]
    where we used $\sum_{m \ge 0} q^m = 1/(1-q)$ and $\sum_{m \ge 0} m q^m = q/(1-q)^2$. This proves \cref{equation:corollary-V_N-convergence-rate-1}.
    By the proof of \Cref{corollary:V-defined}, $|V(\varpi, x_0)|$ and $|V_N(\varpi, x_0)|$ are bounded by a constant that does not depend on $\varpi$, so $V^*(x_0)$ and $V_N^*(x_0)$ are finite.
    Finally, $V^*(x_0) \le V(\varpi, x_0) \le V_N(\varpi, x_0) + \Delta_N(x_0)$ for every $\varpi \in \Pi$, and taking the infimum over $\varpi$ gives $V^*(x_0) \le V_N^*(x_0) + \Delta_N(x_0)$; exchanging the roles of $V$ and $V_N$ gives $V_N^*(x_0) \le V^*(x_0) + \Delta_N(x_0)$.
\end{proof}

We are now ready to show the existence of an optimal policy.
First, the cost function is bounded below. Let $(x, a) \in \text{Gr}(A)$ with $x = (y, k, r, t)$. The first two terms in the square bracket of \cref{equation:cost-function} are nonnegative, because $f \ge 0$, $a \ge 0$ and the weights are nonnegative (\Cref{assumption:standing}). Moreover, $(1 + r - ||a||_1)^\nu \le (1 + \bar{r})^\nu$, because $r \le \bar{r}$. Hence $c(x, a) \ge -e^{-\gamma t}w_R(1+\bar{r})^\nu$.
Rewrite $V(\varpi, x_0)$ as
\begin{equation} \label{equation:V-hat}
    \begin{split}
        &V(\varpi, x_0) = \mathbb{E}\sum_n \left\{[c(X_n, \xi_n) + e^{-\gamma (t_0 + \tau_n)}w_R(1+\bar{r})^\nu] - e^{-\gamma (t_0 + \tau_n)}w_R(1+\bar{r})^\nu \right\} \\
        &= \mathbb{E}\sum_n [c(X_n, \xi_n) + e^{-\gamma (t_0+\tau_n)}w_R(1 + \bar{r})^\nu] - \mathbb{E}\sum_n e^{-\gamma (t_0 + \tau_n)} w_R(1 + \bar{r})^\nu.
    \end{split}
\end{equation}
We could freely split the series and expectation above because of \Cref{corollary:V-defined} and $\mathbb{E}\sum_n e^{-\gamma (t_0 + \tau_n)} \le e^{-\gamma t_0}\sum_n \left(\frac{\bar{\lambda}}{\bar{\lambda} + \gamma} \right)^n < \infty$.
Let
\[
    D(x_0) := w_R(1+\bar{r})^\nu\, \mathbb{E}\sum_{n=0}^\infty e^{-\gamma(t_0+\tau_n)}, \qquad x_0 = (y_0, k_0, r_0, t_0) \in \mathbb{X},
\]
which is finite by the bound just mentioned. The signal times are not controlled: by \cref{equation:signal-construction}, they are built from $k_0$, $(E_n)_{n \ge 1}$ and $(U_n)_{n \ge 1}$ only. Hence $D(x_0)$ does not depend on the policy, and it depends on $x_0$ only through $(k_0, t_0)$; more precisely, $D(y_0, k_0, r_0, t_0) = e^{-\gamma t_0} D(y_0, k_0, r_0, 0)$. Define the modified cost
\begin{equation} \label{equation:c-hat}
    \hat{c}(x, a) := c(x, a) + e^{-\gamma t}w_R(1+\bar{r})^\nu \ge 0, \qquad (x, a) \in \text{Gr}(A),
\end{equation}
and let $\hat{V}(\varpi, x_0) := \mathbb{E}\sum_{n=0}^\infty \hat{c}(X_n, \xi_n)$ and $\hat{V}^*(x_0) := \inf_{\varpi \in \Pi} \hat{V}(\varpi, x_0)$. Then \cref{equation:V-hat} reads $V(\varpi, x_0) = \hat{V}(\varpi, x_0) - D(x_0)$ for every $\varpi \in \Pi$. Hence $V^*(x_0) = \hat{V}^*(x_0) - D(x_0)$, and a policy is optimal for $V$ if and only if it is optimal for $\hat{V}$. This is an important reduction: the original problem is equivalent to the one with the nonnegative cost $\hat{c}$. Moreover, $\hat{c}$ is continuous on $\text{Gr}(A)$, because $f$ is continuous (\Cref{assumption:standing}) and $1 + r - ||a||_1 \ge 1$ on $\text{Gr}(A)$.

\begin{theorem}[Existence of Optimal Policy] \label{theorem:existence-of-optimal-policy}
    Suppose the assumption of \Cref{theorem:one-stage-cost-estimate} holds, and let $\mathbb{A} := \{a \in \mathbb{R}_+^d : ||a||_1 \le \bar{r}\}$. Then:
    \begin{enumerate}[(i)]
        \item The value function $\hat{V}^*$ is lower semicontinuous on $\mathbb{X}$ and satisfies the optimality equation
        \begin{equation} \label{equation:theorem-existence-of-optimal-policy-1}
            \hat{V}^*(x) = \min_{a \in A(x)}\Big\{\hat{c}(x, a) + \int_\mathbb{X} \hat{V}^*(z)\, Q(\dd{z} | x, a)\Big\}, \quad x \in \mathbb{X},
        \end{equation}
        where the minimum is attained.
        \item There exists a stationary deterministic policy $\varpi^* : \mathbb{X} \rightarrow \mathbb{A}$ (\Cref{definition:policies}) such that, for every $x \in \mathbb{X}$, the action $\varpi^*(x)$ attains the minimum in \cref{equation:theorem-existence-of-optimal-policy-1}. Every such policy is optimal:
        \begin{equation} \label{equation:theorem-existence-of-optimal-policy-2}
            \hat{V}(\varpi^*, x) = \hat{V}^*(x), \quad x \in \mathbb{X}.
        \end{equation}
        \item The value function $V^* = \hat{V}^* - D$ of the original problem is lower semicontinuous on $\mathbb{X}$, and every policy $\varpi^*$ as in (ii) is optimal for it: $V(\varpi^*, x) = V^*(x)$ for all $x \in \mathbb{X}$.
    \end{enumerate}
\end{theorem}

\begin{proof}
    We apply Theorem 2 of \textcite{Feinberg} to the Markov decision model with state space $\mathbb{X}$, action space $\mathbb{A}$, feasible action sets $A(x)$, one-step cost $\hat{c}$, transition kernel $Q$ and discount factor $1$. Their results are stated for discount factors in $[0, 1)$. By the remark after the proof of their Theorem 2, the conclusions remain correct for discount factor $1$ when the one-step cost is nonnegative, as $\hat{c}$ is. We first check their setting and their Assumption (W$^*$), and then relate their policies to ours.

    \emph{Setting.} The state space $\mathbb{X}$ is a closed subset of $\mathbb{R}^d \times \mathscr{K} \times \mathbb{R} \times \mathbb{R}$, and $\mathbb{A}$ is a compact subset of $\mathbb{R}^d$, so both are Borel subsets of Polish spaces. Let $x = (y, k, r, t) \in \mathbb{X}$. Then $A(x) \subseteq \mathbb{A}$, because $||a||_1 \le r \le \bar{r}$ for $a \in A(x)$. Also, $A(x)$ is nonempty because $0 \in A(x)$; in particular, $x \mapsto 0$ is a measurable selector of $A$. The graph $\text{Gr}(A)$ is closed (\Cref{definition:feasible-action-set}), hence Borel. The cost $\hat{c}$ is defined on $\text{Gr}(A)$, as \textcite{Feinberg} require. Their transition kernel, however, must be defined on all of $\mathbb{X} \times \mathbb{A}$, while $Q$ is defined on $\text{Gr}(A)$ (\Cref{theorem:transition-kernel}). We extend it by $Q(\cdot | x, a) := Q(\cdot | x, 0)$ for $a \notin A(x)$, which is possible because $(x, 0) \in \text{Gr}(A)$. The extension is a stochastic kernel on $\mathbb{X}$ given $\mathbb{X} \times \mathbb{A}$. It does not change $Q$ on $\text{Gr}(A)$, and only the values on $\text{Gr}(A)$ enter their assumptions, their results and the dynamics of the process.

    \emph{Assumption (W$^*$).} We check its three conditions.
    \begin{enumerate}[(a)]
        \item $\hat{c}$ is lower semicontinuous and bounded below on $\text{Gr}(A)$: it is continuous and nonnegative there (see \cref{equation:c-hat} and the sentence after it).
        \item If $x_m \rightarrow x$ in $\mathbb{X}$ and $a_m \in A(x_m)$ for every $m$, then $(a_m)$ has a limit point in $A(x)$. (The condition of \textcite{Feinberg} only concerns sequences along which $\hat{c}(x_m, a_m)$ is bounded above, but we do not need this.) Indeed, write $x_m = (y_m, k_m, r_m, t_m)$ and $x = (y, k, r, t)$. Since $a_m \in \mathbb{A}$ and $\mathbb{A}$ is compact, a subsequence $(a_{m'})$ converges to some $a \in \mathbb{A}$; in particular, $a \in \mathbb{R}_+^d$. Since $||a_{m'}||_1 \le r_{m'}$ and $r_{m'} \rightarrow r$, we get $||a||_1 \le r$. Since $\mathscr{K}$ is discrete, $k_m = k$ for all large $m$. Hence $a_{m',j} = 0$ for $j \notin k$ and all large $m'$, and so $a_j = 0$ for $j \notin k$. Thus $a \in A(x)$.
        \item $Q$ is weakly continuous on $\text{Gr}(A)$: this is \Cref{theorem:weak-continuity}.
    \end{enumerate}

    \emph{Conclusions of \textcite{Feinberg}.} For $x \in \mathbb{X}$, let $v(x)$ be the infimum, over their policies, of their expected total cost $\mathbb{E}\sum_{n \ge 0} \hat{c}(x_n, a_n)$ from the initial state $x$. By Theorem 2(i), (iv) and (v) of \textcite{Feinberg}:
    \begin{itemize}
        \item $v$ is lower semicontinuous and satisfies \cref{equation:theorem-existence-of-optimal-policy-1} with $v$ in place of $\hat{V}^*$, with the minimum attained;
        \item there is a measurable function $\varpi^*: \mathbb{X} \rightarrow \mathbb{A}$ such that, for every $x$, $\varpi^*(x) \in A(x)$ and $\varpi^*(x)$ attains this minimum;
        \item every such stationary policy is optimal for their problem.
    \end{itemize}

    \emph{Relation between the two problems.} Fix $\varpi \in \Pi$ and $x_0 \in \mathbb{X}$. By \Cref{remark:policies}, the stochastic kernels $h_n \mapsto \P\{\varpi_n(h_n, \chi_n) \in \cdot\}$ form a policy of \textcite{Feinberg}, and they give the conditional distribution of $\xi_n$ given $H_n$. Given $(H_n, \xi_n)$, the next state $X_{n+1}$ has the distribution $Q(\cdot | X_n, \xi_n)$ (\Cref{theorem:markov-property,theorem:transition-kernel}). By the Ionescu-Tulcea theorem (see, e.g., \textcite[Proposition~C.10 and Remark~C.11]{Lerma1}), these conditional distributions and $X_0 = x_0$ determine the law of the sequence $(X_0, \xi_0, X_1, \xi_1, \dots)$. This sequence therefore has the same law as the state-action sequence of \textcite{Feinberg} under the induced policy with initial state $x_0$. Since $\hat{c} \ge 0$, the two expected total costs are equal, so $\hat{V}(\varpi, x_0) \ge v(x_0)$. Taking the infimum over $\varpi \in \Pi$ gives $\hat{V}^* \ge v$. Conversely, a stationary deterministic policy $\varpi^*$ as in the previous paragraph belongs to $\Pi$ (\Cref{definition:policies}). The same argument gives $\hat{V}(\varpi^*, x) = v(x)$, so $\hat{V}^* \le v$. Hence $\hat{V}^* = v$, which proves (i) and (ii).

    Finally, $D(x)$ depends on $x = (y, k, r, t)$ only through $k$, which ranges over a discrete set, and through the factor $e^{-\gamma t}$, which is continuous in $t$. Hence $D$ is continuous on $\mathbb{X}$, and $V^* = \hat{V}^* - D$ is lower semicontinuous. For every policy $\varpi^*$ as in (ii), $V(\varpi^*, x) = \hat{V}(\varpi^*, x) - D(x) = \hat{V}^*(x) - D(x) = V^*(x)$, which proves (iii).
\end{proof}

\begin{remark} \label{remark:time-independence}
    An optimal injection rule does not need to depend on calendar time. For $x = (y, k, r, t) \in \mathbb{X}$, write $x^0 := (y, k, r, 0)$. Then
    \begin{equation} \label{equation:remark-time-independence-1}
        \hat{V}^*(x) = e^{-\gamma t}\hat{V}^*(x^0) \quad \text{and} \quad V^*(x) = e^{-\gamma t}V^*(x^0),
    \end{equation}
    and the actions that attain the minimum in \cref{equation:theorem-existence-of-optimal-policy-1} at $x$ and at $x^0$ are the same. Consequently, if $\varpi^*$ is as in \Cref{theorem:existence-of-optimal-policy}(ii), then so is the policy $x \mapsto \varpi^*(x^0)$. This policy is optimal and depends only on the members' wealth, the set of members directly hit by the most recent signal and the fund's balance.

    To see this, note that $\hat{c}(x, a) = e^{-\gamma t}\hat{c}(x^0, a)$ by \cref{equation:cost-function,equation:c-hat}, and $A(x) = A(x^0)$. Moreover, in \cref{equation:theorem-transition-kernel-4} only the time component of $Z^{(x, a)}$ depends on $t$, and it equals $t + S$. So, by \cref{equation:theorem-transition-kernel-5}, $Q(\cdot | x, a)$ is the image of $Q(\cdot | x^0, a)$ under the shift $(y', k', r', t') \mapsto (y', k', r', t' + t)$. By Theorem 2(i) and (ii) of \textcite{Feinberg} (with discount factor $1$, as in the proof of \Cref{theorem:existence-of-optimal-policy}), the functions defined by $v_0 := 0$ and
    \[
        v_{n+1}(x) := \min_{a \in A(x)}\Big\{\hat{c}(x, a) + \int_\mathbb{X} v_n(z)\, Q(\dd{z} | x, a)\Big\}, \qquad n \ge 0,
    \]
    converge to $\hat{V}^*$. Suppose that $v_n(x) = e^{-\gamma t}v_n(x^0)$ for all $x \in \mathbb{X}$, which holds for $n = 0$. Then $v_n(y', k', r', t' + t) = e^{-\gamma t}v_n(z')$ for every $z' = (y', k', r', t') \in \mathbb{X}$, so
    \[
        \int_\mathbb{X} v_n(z)\, Q(\dd{z} | x, a) = \int_\mathbb{X} v_n(y', k', r', t' + t)\, Q(\dd{z'} | x^0, a) = e^{-\gamma t}\int_\mathbb{X} v_n(z')\, Q(\dd{z'} | x^0, a),
    \]
    and hence $v_{n+1}(x) = e^{-\gamma t}v_{n+1}(x^0)$. By induction and letting $n \rightarrow \infty$, we obtain the first identity in \cref{equation:remark-time-independence-1}. The second follows from $V^* = \hat{V}^* - D$ and $D(x) = e^{-\gamma t}D(x^0)$. Finally, the same computation with $\hat{V}^*$ in place of $v_n$ shows that the expression minimized in \cref{equation:theorem-existence-of-optimal-policy-1} at $x$ is $e^{-\gamma t}$ times the one at $x^0$. Both are minimized over the same set $A(x) = A(x^0)$, so they have the same minimizers.
\end{remark}

In terms of the application, calendar time enters the problem only through the discount factor $e^{-\gamma t}$: the shock rates, the distribution of the losses, the dynamics of the wealth and of the fund, and the cost are the same at every date. The fund's manager therefore faces the same decision problem at every market signal. The injection depends on the situation of the group when the shock arrives (how much each member has, who was hit and how much is left in the fund), but not on when the shock arrives. In particular, a microfinance institution can compute the injection rule once, as a function of $(y, k, r)$, and apply it at every shock, and computing it numerically involves one state variable fewer. Similarly, $V^*(x^0)$ is the expected total cost from the current signal onwards, measured in money of the current date, and the factor $e^{-\gamma t}$ only converts it into money of time $0$. This relies on the model being time-homogeneous: the shock rates and the other parameters do not depend on calendar time, and costs are discounted at a constant rate $\gamma$. If shocks were seasonal, for instance if droughts were more likely in the dry season, the rates $\lambda_{k,k'}$ would depend on the time of year, and so would the optimal injections; we leave this extension to future work.

\printbibliography

@article{JassoFuentesVasquez2025,
    author = {Jasso-Fuentes, H{\'e}ctor and Vasquez-Rojas, Fidel},
    title = {Optimal stopping problems with constraints for continuous-time general {M}arkov processes},
    journal = {ESAIM: Control, Optimisation and Calculus of Variations},
    year = 2025,
    volume = 31,
    eid = {53},
    doi = {10.1051/cocv/2025035},
}

@book{Paolo,
    author = {Baldi, Paolo},
    title = {Stochastic Calculus: An Introduction Through Theory and Exercises},
    publisher = {Springer},
    year = 2017,
    doi = {10.1007/978-3-319-62226-2}
}

@article{Menaldi1,
    author = {Bayraktar, Erhan and Emmerling, Thomas and Menaldi, José-Luis},
    title = {On the Impulse Control of Jump Diffusions},
    journal = {SIAM Journal on Control and Optimization},
    year = 2013,
    volume = {51},
    number = {3},
    pages = {2612--2637},
    doi = {10.1137/120863836}
}

@article{Kallenberg,
    author = {Kallenberg, Olav},
    title = {On the Existence of Universal Functional Solutions to Classical SDE's},
    journal = {The Annals of Probability},
    year = 1996,
    volume = {24},
    number = {1},
    pages = {196--205},
    doi = {10.1214/aop/1042644713}
}

@book{Mao,
    author = {Mao, Xuerong and Yuan, Chenggui},
    title = {Stochastic Differential Equations with Markovian Switching},
    publisher = {Imperial College Press},
    year = 2006,
    doi = {10.1142/p473}
}

@book{Lerma1,
    author = {Hernández-Lerma, Onésimo and Lasserre, Jean-Bernard},
    title = {Discrete-Time Markov Control Processes - Basic Optimality Criteria},
    publisher = {Springer},
    year = 1996,
    doi = {10.1007/978-1-4612-0729-0}
}

@book{Skorohod,
    author = {Skorokhod, Anatolii Vladimirovich},
    title = {Asymptotic Methods in the Theory of Stochastic Differential Equations},
    series = {Translations of Mathematical Monographs},
    volume = {78},
    publisher = {American Mathematical Society},
    year = {1989},
    doi = {10.1090/mmono/078}
}

@article{Ghosh,
    author = {Ghosh, Mrinal K. and Arapostathis, Aristotle and Marcus, Steven I.},
    title = {Ergodic Control of Switching Diffusions},
    journal = {SIAM Journal on Control and Optimization},
    year = {1997},
    volume = {35},
    number = {6},
    pages = {1952--1988},
    doi = {10.1137/S0363012996299302}
}

@book{Ikeda,
    author = {Ikeda, Nobuyuki and Watanabe, Shinzo},
    title = {Stochastic Differential Equations and Diffusion Processes},
    edition = {2},
    series = {North-Holland Mathematical Library},
    volume = {24},
    publisher = {North-Holland Publishing Company and Kodansha},
    year = {1989},
    isbn = {0-444-87378-3}
}

@article{Feinberg,
    author = {Feinberg, Eugene A. and Kasyanov, Pavlo O. and Zadoianchuk, Nina V.},
    title = {Average Cost Markov Decision Processes with Weakly Continuous Transition Probabilities},
    journal = {Mathematics of Operations Research},
    year = {2012},
    volume = {37},
    number = {4},
    pages = {591--607},
    doi = {10.1287/moor.1120.0555}
}

@article{ProtterQuintos2022,
    author = {Protter, Philip and Quintos, Alejandra},
    title = {Optimal group size in microlending},
    journal = {Annals of Finance},
    year = {2022},
    volume = {18},
    pages = {121--132},
    doi = {10.1007/s10436-020-00382-0}
}

@article{JassoFuentesQuintosYang2026,
    author = {Jasso-Fuentes, Héctor and Quintos, Alejandra and Yang, Xinta},
    title = {Group survival probability under contagion in microlending},
    journal = {Annals of Finance},
    year = {2026},
    volume = {22},
    eid = {7},
    doi = {10.1007/s10436-026-00482-3}
}

@article{FosterGreerThorbecke1984,
  author  = {Foster, James and Greer, Joel and Thorbecke, Erik},
  title   = {A Class of Decomposable Poverty Measures},
  journal = {Econometrica},
  volume  = {52},
  number  = {3},
  pages   = {761--766},
  year    = {1984},
  doi     = {10.2307/1913475}
}

@book{Norris,
    author = {Norris, James R.},
    title = {Markov Chains},
    series = {Cambridge Series in Statistical and Probabilistic Mathematics},
    publisher = {Cambridge University Press},
    year = {1997},
    doi = {10.1017/CBO9780511810633},
    location = {Cambridge}
}

@book{Yunus2007,
    author = {Yunus, Muhammad},
    note = {With Alan Jolis},
    title = {Banker to the Poor: Micro-Lending and the Battle Against World Poverty},
    publisher = {PublicAffairs},
    address = {New York},
    year = {2007},
    isbn = {978-1-58648-198-8}
}

@article{DupuisWang2002,
    author = {Dupuis, Paul and Wang, Hui},
    title = {Optimal Stopping with Random Intervention Times},
    journal = {Advances in Applied Probability},
    year = {2002},
    volume = {34},
    pages = {141--157},
    doi = {10.1239/aap/1019160954}
}

@article{MenaldiRobin2016,
    author = {Menaldi, Jos{\'e}-Luis and Robin, Maurice},
    title = {On Some Optimal Stopping Problems with Constraint},
    journal = {SIAM Journal on Control and Optimization},
    year = {2016},
    volume = {54},
    number = {5},
    pages = {2650--2671},
    doi = {10.1137/15M1040001}
}

@article{MenaldiRobin2017,
    author = {Menaldi, Jos{\'e}-Luis and Robin, Maurice},
    title = {On Some Impulse Control Problems with Constraint},
    journal = {SIAM Journal on Control and Optimization},
    year = {2017},
    volume = {55},
    number = {5},
    pages = {3204--3225},
    doi = {10.1137/16M1090302}
}

@book{Davis1993,
    author = {Davis, M. H. A.},
    title = {Markov Models and Optimization},
    publisher = {Chapman \& Hall},
    address = {London},
    year = {1993},
    isbn = {978-0-412-31410-0}
}

@article{AlbrecherCheungThonhauser2011,
    author = {Albrecher, Hansj{\"o}rg and Cheung, Eric C. K. and Thonhauser, Stefan},
    title = {Randomized Observation Periods for the Compound {Poisson} Risk Model: Dividends},
    journal = {ASTIN Bulletin},
    year = {2011},
    volume = {41},
    number = {2},
    pages = {645--672},
    doi = {10.2143/AST.41.2.2136991}
}

@article{AvramPerezYamazaki2018,
    author = {Avram, Florin and P{\'e}rez, Jos{\'e}-Luis and Yamazaki, Kazutoshi},
    title = {Spectrally Negative {L\'evy} Processes with {Parisian} Reflection Below and Classical Reflection Above},
    journal = {Stochastic Processes and their Applications},
    year = {2018},
    volume = {128},
    number = {1},
    pages = {255--290},
    doi = {10.1016/j.spa.2017.04.013}
}

@article{NobaPerezYamazakiYano2018,
    author = {Noba, Kei and P{\'e}rez, Jos{\'e}-Luis and Yamazaki, Kazutoshi and Yano, Kouji},
    title = {On Optimal Periodic Dividend and Capital Injection Strategies for Spectrally Negative {L\'evy} Models},
    journal = {Journal of Applied Probability},
    year = {2018},
    volume = {55},
    pages = {1272--1286},
    doi = {10.1017/jpr.2018.85}
}

@misc{MataMorenoFrancoNobaPerez2022,
    author = {Mata, Dante and Moreno-Franco, Harold A. and Noba, Kei and P{\'e}rez, Jos{\'e}-Luis},
    title = {On the Bailout Dividend Problem with Periodic Dividend Payments for Spectrally Negative {Markov} Additive Processes},
    year = {2022},
    eprint = {2207.01126},
    eprinttype = {arXiv},
    eprintclass = {math.PR}
}

@article{FouqueIchiba2013,
    author = {Fouque, Jean-Pierre and Ichiba, Tomoyuki},
    title = {Stability in a Model of Interbank Lending},
    journal = {SIAM Journal on Financial Mathematics},
    year = {2013},
    volume = {4},
    pages = {784--803},
    doi = {10.1137/110841096}
}

@article{CarmonaFouqueSun2015,
    author = {Carmona, Ren{\'e} and Fouque, Jean-Pierre and Sun, Li-Hsien},
    title = {Mean Field Games and Systemic Risk},
    journal = {Communications in Mathematical Sciences},
    year = {2015},
    volume = {13},
    number = {4},
    pages = {911--933},
    doi = {10.4310/CMS.2015.v13.n4.a4}
}

@inproceedings{DienerDienerKhodrProtter2009,
    author = {Diener, Francine and Diener, Marc and Khodr, Osman and Protter, Philip},
    title = {Mathematical Models for Microlending},
    booktitle = {Proceedings of the 16th Mathematical Conference of the Bangladesh Mathematical Society},
    address = {Dhaka, Bangladesh},
    year = {2009}
}

@article{JarrowProtter2019,
    author = {Jarrow, Robert and Protter, Philip},
    title = {Fair Microfinance Loan Rates},
    journal = {International Review of Finance},
    year = {2019},
    volume = {19},
    number = {4},
    pages = {909--918},
    doi = {10.1111/irfi.12195}
}

@article{BoiquayeProtter2024,
    author = {Boiquaye, Perpetual Andam and Protter, Philip},
    title = {Probability of No Default for a Microloan under Uncertainty},
    journal = {Annals of Finance},
    year = {2024},
    volume = {20},
    pages = {521--528},
    doi = {10.1007/s10436-024-00455-4}
}

@article{BoiquayeKorvehGyamerah2026,
    author = {Boiquaye, Perpetual Andam and Korveh, Edward and Gyamerah, Samuel Asante},
    title = {A {Markov} Process Model of Joint Liability and Loan Repayment for Sustainable Microfinance},
    journal = {Annals of Finance},
    year = {2026},
    volume = {22},
    eid = {1},
    doi = {10.1007/s10436-026-00478-z}
}

@article{AlbrecherBauerleThonhauser2011,
    author = {Albrecher, Hansj{\"o}rg and B{\"a}uerle, Nicole and Thonhauser, Stefan},
    title = {Optimal Dividend-Payout in Random Discrete Time},
    journal = {Statistics \& Risk Modeling},
    year = {2011},
    volume = {28},
    number = {3},
    pages = {251--276},
    doi = {10.1524/stnd.2011.1097}
}

\begin{appendices}
    \section{Miscellaneous Results}
    Here we collect some useful results that are used as lemmas in the main body of this paper.
    For results that are well known, we only sketch the proofs.

    \begin{lemma}[Freezing Lemma] \label{lemma:freezing-lemma}
        Let $(\Omega, \mathcal{F}, \P)$ be a probability space, and let $\mathcal{H}$ and $\mathcal{D}$ be independent sub-$\sigma$-fields of $\mathcal{F}$.
        Let $X$ be a $\mathcal{D}$-measurable random variable with values in a measurable space $(E, \mathcal{E})$, and let $g: E \times \Omega \rightarrow \mathbb{R}$ be $\mathcal{E} \otimes \mathcal{H}/\mathscr{B}(\mathbb{R})$-measurable and either bounded or nonnegative.
        Then
        \begin{equation} \label{equation:lemma-freezing-lemma-1}
            \mathbb{E}\{g(X, \cdot) | \mathcal{D}\} = \bar{g}(X), \qquad \text{where } \bar{g}(x) := \mathbb{E}\{g(x, \cdot)\}, \ x \in E.
        \end{equation}
        If $g$ is nonnegative, no integrability is needed, and both sides of \cref{equation:lemma-freezing-lemma-1} take values in $[0, \infty]$.
    \end{lemma}
    \begin{proof}
        For bounded $g$, this follows from Lemma 4.1 of \textcite{Paolo}. If $g$ is nonnegative, apply the bounded case to $g \wedge m$, $m \ge 1$, and let $m \rightarrow \infty$: the left-hand side converges to $\mathbb{E}\{g(X, \cdot) | \mathcal{D}\}$ by the monotone convergence theorem for conditional expectations, and $\mathbb{E}\{(g \wedge m)(x, \cdot)\} \uparrow \bar{g}(x)$ for every $x \in E$ by monotone convergence.
    \end{proof}

    \begin{lemma}[Measurability of Evaluation Map] \label{lemma:measurability-of-evaluation-map}
        For any dimension $d \in \mathbb{Z}_+$, define the evaluation map $\mathfrak{e}: \mathbb{R}_+ \times \mathscr{C}(\mathbb{R}_+, \mathbb{R}^d) \rightarrow \mathbb{R}^d$ via $\mathfrak{e}(t, y) := y(t)$.
        Claim function $\mathfrak{e}$ is $\mathscr{B}(\mathbb{R}_+) \otimes \mathscr{B}(\mathscr{C}(\mathbb{R}_+, \mathbb{R}^d))/\mathscr{B}(\mathbb{R}^d)$-measurable.
    \end{lemma}
    \begin{proof}[\textit{Proof Sketch}]
        We first show $\mathfrak{e}$ is continuous with respect to the metric topology on $\mathbb{R}_+ \times \mathscr{C}(\mathbb{R}_+, \mathbb{R}^d)$.
        Take $(t_n, y_n) \rightarrow (t, y)$.
        Because $y_n$ converges to $y$ uniformly on compacts (refer to the standard metrization of $\mathscr{C}(\mathbb{R}_+, \mathbb{R}^d)$), one could show $y_n(t_n) \rightarrow y(t)$.
        Next, because $\mathbb{R}_+$ and $\mathscr{C}(\mathbb{R}_+, \mathbb{R}^d)$ are Polish spaces, $\mathscr{B}(\mathbb{R}_+) \otimes \mathscr{B}(\mathscr{C}(\mathbb{R}_+, \mathbb{R}^d)) = \mathscr{B}(\mathbb{R}_+ \times \mathscr{C}(\mathbb{R}_+, \mathbb{R}^d))$.
        Combined with the continuity property of $\mathfrak{e}$, the claim follows.
    \end{proof}

    \begin{lemma} \label{lemma:markov-kappa}
        The chain $(\kappa_{\tau_n})$ satisfies Markov property, i.e.
        \begin{equation} \label{equation:lemma-markov-kappa-1}
            \P\{\kappa_{\tau_{n+1}} = k' | H_n\} = \P\{\kappa_{\tau_{n+1}} = k' | \kappa_{\tau_n}\}.
        \end{equation}
    \end{lemma}
    \begin{proof}[\textit{Proof Sketch}]

        By \cref{equation:signal-construction}, $\kappa_{\tau_{n+1}} = F(\kappa_{\tau_n}, U_{n+1})$ and $\tau_{n+1} - \tau_n = E_{n+1}/\lambda_{\kappa_{\tau_n}}$.
        The pair $(E_{n+1}, U_{n+1})$ is independent of $(H_n, \xi_n)$: by \Cref{definition:policies}, $H_n$ and $\xi_n$ are functions of $B$, $(E_i, U_i, \eta_i)_{i \le n}$ and $(\chi_i)_{i \le n}$, which are independent of $(E_{n+1}, U_{n+1})$ by \Cref{assumption:standing}.
        Since $\kappa_{\tau_n}$ is a component of $X_n$, \nCref{lemma:freezing-lemma} gives
        \[
            \P\{\kappa_{\tau_{n+1}} = k' | H_n\} = \frac{\lambda_{\kappa_{\tau_n}, k'}}{\lambda_{\kappa_{\tau_n}}},
        \]
        because $\P\{F(k, U_{n+1}) = k'\} = \lambda_{k,k'}/\lambda_k$ for each $k \in \mathscr{K}$.
        The right-hand side is a function of $\kappa_{\tau_n}$, so by the tower property it is also equal to $\P\{\kappa_{\tau_{n+1}} = k' | \kappa_{\tau_n}\}$.
    \end{proof}

    \begin{lemma} \label{lemma:signal-interval-distribution}
        The conditional distribution of market signals' inter-arrival times is given by
        \begin{equation} \label{lemma-signal-interval-distribution-1}
            \P\{\tau_{n+1}-\tau_n > s | H_n\} = e^{-s\lambda_{\kappa_{\tau_n}}}.
        \end{equation}
    \end{lemma}
    \begin{proof}[\textit{Proof Sketch}]

        By \cref{equation:signal-construction}, $\tau_{n+1} - \tau_n = E_{n+1}/\lambda_{\kappa_{\tau_n}}$, where $E_{n+1}$ is exponential with mean $1$ and independent of $H_n$ (see the proof of \Cref{lemma:markov-kappa}).
        By \nCref{lemma:freezing-lemma}, $\P\{\tau_{n+1}-\tau_n > s | H_n\} = e^{-s\lambda_{\kappa_{\tau_n}}}$, since $\P\{E_{n+1} > s\lambda_k\} = e^{-s\lambda_k}$ for each $k \in \mathscr{K}$.
    \end{proof}

    \begin{remark} \label{remark:markov-kappa}
        The same argument used in \Cref{lemma:markov-kappa} also produces
        \begin{equation} \label{equation:remark-markov-kappa-1}
            \mathbb{E}\{1_\Gamma(\kappa_{\tau_{n+1}}, \tau_{n+1}-\tau_n) | H_n\} = \mathbb{E}\{1_\Gamma(\kappa_{\tau_{n+1}}, \tau_{n+1}-\tau_n) | \kappa_{\tau_n}\},
        \end{equation}
        for Borel sets $\Gamma \subset \mathscr{K}\times \mathbb{R}_+$.
    \end{remark}

    \begin{lemma} \label{lemma:more-on-signal-interval-distribution}
        If we also specify the mark $\kappa_{\tau_{n+1}}$ of the next signal, the above result can be strengthened to
        \begin{equation} \label{equation:lemma-more-on-signal-interval-distribution-1}
            \P\{\tau_{n+1}-\tau_n \le s, \kappa_{\tau_{n+1}}=k' | H_n\} = \frac{\lambda_{\kappa_{\tau_n}, k'}}{\lambda_{\kappa_{\tau_n}}}(1 - e^{-s\lambda_{\kappa_{\tau_n}}}).
        \end{equation}
    \end{lemma}
    \begin{proof}[\textit{Proof Sketch}]

        By \cref{equation:signal-construction}, $\tau_{n+1} - \tau_n = E_{n+1}/\lambda_{\kappa_{\tau_n}}$ and $\kappa_{\tau_{n+1}} = F(\kappa_{\tau_n}, U_{n+1})$, where $E_{n+1}$, $U_{n+1}$ and $H_n$ are independent (see the proof of \Cref{lemma:markov-kappa}).
        For fixed $k \in \mathscr{K}$, the set $\{u \in [0, 1) : F(k, u) = k'\}$ is an interval of length $\lambda_{k,k'}/\lambda_k$ (empty if $\lambda_{k,k'} = 0$), so
        \[
            \P\{E_{n+1}/\lambda_k \le s, F(k, U_{n+1}) = k'\} = \frac{\lambda_{k,k'}}{\lambda_k}(1 - e^{-s\lambda_k}).
        \]
        The result follows from \nCref{lemma:freezing-lemma} with $X = \kappa_{\tau_n}$.
    \end{proof}

    \begin{remark} \label{remark:more-on-signal-interval-distribution}
        It is also useful to consider the following variation of the above lemma:
        \begin{equation} \label{equation:lemma-more-on-signal-interval-distribution-4}
            \P\{\tau_{n+1}-\tau_n \in \dd{s}, \kappa_{\tau_{n+1}} = k' | H_n\} = \lambda_{\kappa_{\tau_n},k'}e^{-s\lambda_{\kappa_{\tau_n}}} \dd{s}.
        \end{equation}
    \end{remark}

    \begin{remark} \label{remark:further-on-signal-interval-distribution}
        If we condition an extra $\xi_n$ on the left hand side of \cref{equation:lemma-more-on-signal-interval-distribution-4}, the right hand side stays the same, because $(E_{n+1}, U_{n+1})$ is independent of $(H_n, \xi_n)$ (see the proof of \Cref{lemma:markov-kappa}), so the proof of \Cref{lemma:more-on-signal-interval-distribution} applies with $(H_n, \xi_n)$ in place of $H_n$.
    \end{remark}

    \begin{lemma} \label{lemma:SDE-convergence}
        Given the usual SDE $\dd{Y}_t = b(Y_t) \dd{t} + \sigma(Y_t) \dd{B_t}$, where $b, \sigma$ satisfy linear growth and Lipschitz conditions.
        Denote by $Y^y$ and $Y^{y_n}$ the SDE solution with initial condition $y$ and $y_n$ respectively, both driven by the same Brownian motion $B$.
        Then
        \begin{equation} \label{equation:lemma-SDE-convergence}
            \mathbb{E}\left\{\sup_{0\le s \le t} ||Y^{y_n}_s - Y^y_s||^2\right\} \le C||y_n-y||^2,
        \end{equation}
        where the constant $C$ depends only on $t$ and the Lipschitz constants of $b, \sigma$.
    \end{lemma}
    \begin{proof}[\textit{Proof Sketch}]
        Follow the standard arguments of Burkholder-Davis-Gundy's Inequality and Gronwall's Inequality.
        These arguments, applied to the difference $Y^{y_n} - Y^y$, use only the Lipschitz conditions on $b$ and $\sigma$; the linear growth condition only guarantees that the moments involved are finite. Hence $C$ depends neither on $y$ nor on $y_n$. A related estimate, for SDEs with Markovian switching whose solutions may start at different times and from bounded initial points, is \textcite[Lemma~3.3]{Mao}; there the constant also depends on the bound on the initial points.
    \end{proof}

    \begin{lemma} \label{lemma:SDE-estimate}
        Consider the usual SDE $\dd{Y_t} = b(Y_t) \dd{t} + \sigma(Y_t) \dd{B_t}$ with $Y_t, B_t \in \mathbb{R}^d, \sigma:\mathbb{R}^d \rightarrow \mathbb{R}^{d \times d}$, and $b, \sigma$ satisfying the usual Lipschitz and linear growth condition.
        The initial condition $Y_0 \in \mathbb{R}^d$ is deterministic.
        Suppose there exists a constant $\zeta \ge 0$ such that
        \begin{equation} \label{equation:lemma-SDE-estimate-1}
            y^\intercal b(y) + \frac{1}{2}||\sigma(y)||^2 \le \zeta(1 + ||y||^2),
        \end{equation}
        where $||\sigma(y)||$ is the Frobenius norm of the matrix $\sigma(y)$.
        Then, for every $t \ge 0$,
        \begin{equation} \label{equation:lemma-SDE-estimate-2}
            \mathbb{E}\{1 + ||Y_t||^2\} \le (1 + ||Y_0||^2)e^{2\zeta t}.
        \end{equation}
    \end{lemma}
    \begin{proof}
        By Itô's formula,
        \[
            \dd{(1+||Y_t||^2)} = \big(2Y_t^\intercal b(Y_t) + ||\sigma(Y_t)||^2\big) \dd{t} + 2Y_t^\intercal \sigma(Y_t) \dd{B_t}.
        \]
        For $m > ||Y_0||$, let $T_m := \inf\{t \ge 0 : ||Y_t|| \ge m\}$. Since $||Y_s|| \le m$ for $s \le T_m$ and $\sigma$ is continuous, the stopped stochastic integral $\int_0^{t \wedge T_m} 2Y_s^\intercal \sigma(Y_s) \dd{B_s}$ is a martingale with mean zero. Taking expectations and using \cref{equation:lemma-SDE-estimate-1}, we get
        \[
        \begin{split}
            u_m(t) := \mathbb{E}\{1 + ||Y_{t \wedge T_m}||^2\} &\le 1 + ||Y_0||^2 + 2\zeta\, \mathbb{E}\int_0^{t \wedge T_m} (1 + ||Y_s||^2) \dd{s} \\
            &\le 1 + ||Y_0||^2 + 2\zeta \int_0^t u_m(s) \dd{s}.
        \end{split}
        \]
        For the last step: $Y_s = Y_{s \wedge T_m}$ for $s \le t \wedge T_m$, $\zeta \ge 0$ and the integrand is nonnegative, so the integral over $[0, t \wedge T_m]$ is at most the integral of $1 + ||Y_{s \wedge T_m}||^2$ over $[0, t]$; then we use Fubini's theorem. The function $u_m$ is bounded by $1 + m^2$, and continuous by dominated convergence since $Y$ has continuous paths, so Gronwall's inequality gives $u_m(t) \le (1 + ||Y_0||^2)e^{2\zeta t}$ for every $m$. Finally, $Y$ has continuous paths on $[0, \infty)$, so $T_m \rightarrow \infty$ a.s. and $Y_{t \wedge T_m} \rightarrow Y_t$ a.s. By Fatou's lemma, $\mathbb{E}\{1 + ||Y_t||^2\} \le \liminf_{m \rightarrow \infty} u_m(t) \le (1 + ||Y_0||^2)e^{2\zeta t}$.
    \end{proof}

\end{appendices}

\end{document}